\documentclass[journal]{IEEEtran}
\usepackage{graphics} 
\usepackage{epsfig} 
\usepackage{times} 
\usepackage{amsmath} 
\usepackage{amssymb}  
\usepackage{cite}
\usepackage{epstopdf}
\usepackage{subfigure,color,balance}
\usepackage{verbatim}
\usepackage{cases}
\usepackage{enumerate}
\usepackage{balance}

\usepackage[colorlinks=true,      
			linkcolor=black,      
			citecolor=black,      
			filecolor=black,      
			urlcolor=blue]{hyperref}

\newtheorem{theorem}{Theorem}

\newtheorem{remark}{Remark}
\newtheorem{lemma}{Lemma}
\newtheorem{corollary}{Corollary}

\def\B{{\bf B}}

\def\D{{\bf D}}

\def\L{{\bf L}}
\def\M{{\bf M}}

\def\R{{\bf R}}

\def\0{{\bf 0}}
\def\1{{\bf 1}}

\def\eg{{\em e.g.}}
\def\ie{{\em i.e.}}

\def\tr{\mathsf{T}}

\begin{document}
%
\title{
Smoothing Traffic Flow via Control of \\ Autonomous Vehicles}

\author{Yang Zheng,~\IEEEmembership{Member,~IEEE}, Jiawei Wang,~\IEEEmembership{Student Member,~IEEE}, and Keqiang Li
\thanks{Copyright (c) 20xx IEEE. Personal use of this material is permitted. However, permission to use this material for any other purposes must be obtained from the IEEE by sending a request to pubs-permissions@ieee.org.}
\thanks{This work is supported by National Key R\&D Program of China with 2016YFB0100906. The authors also acknowledge the support from TOYOTA.}
\thanks{Yang Zheng was with the Department of Engineering Science, University of Oxford, Oxford, OX1 2JD, U.K. He is now with the School of Engineering and Applied Sciences, and the Harvard Center for Green Buildings and Cities, Harvard University, Cambridge, MA 02138 (zhengy@g.harvard.edu).}%
\thanks{J.~Wang and K.~Li are with the School of Vehicle and Mobility, Tsinghua University, Beijing, China, and with the Center for Intelligent Connected Vehicles \& Transportation, Tsinghua University, Beijing China. (e-mail: wang-jw18@mails.tsinghua.edu.cn, likq@tsinghua.edu.cn). Corresponding author: Keqiang Li.}%
}

\maketitle

\begin{abstract}
The emergence of autonomous vehicles is expected to revolutionize road transportation in the near future. Although large-scale numerical simulations and small-scale experiments have shown promising results, a comprehensive theoretical understanding to smooth traffic flow via autonomous vehicles is lacking. In this paper, from a control-theoretic perspective, we establish analytical results on the controllability, stabilizability, and reachability of a mixed traffic system consisting of human-driven vehicles and autonomous vehicles in a ring road. We show that the mixed traffic system is not completely controllable, but is stabilizable, indicating that autonomous vehicles can not only suppress unstable traffic waves but also guide the traffic flow to a higher speed. Accordingly, we establish the maximum traffic speed achievable via controlling autonomous vehicles. Numerical results show that the traffic speed can be increased by over 6\% when there are only 5\% autonomous vehicles. We also design an optimal control strategy for autonomous vehicles to actively dampen undesirable perturbations. These theoretical findings validate the high potential of autonomous vehicles to smooth traffic flow.
\end{abstract}

\begin{IEEEkeywords}
Autonomous vehicle, mixed traffic flow, controllability, stabilizability.
\end{IEEEkeywords}

\IEEEpeerreviewmaketitle

\section{Introduction}

\IEEEPARstart{M}{odern} societies are increasingly relying on complex road transportation systems to support our daily mobility needs. In some big cities, the traffic demand is placing a heavy burden on existing transportation infrastructures, sometimes leading to severely congested road networks~\cite{batty2008size}. Traffic congestion not only results in the loss of fuel economy and travel efficiency, but also increases the potential risk of traffic accidents and public health~\cite{levy2010evaluation}.

Understanding traffic dynamics is essential if we are to redesign infrastructures, or to guide/control transportation, to mitigate road congestions and smooth traffic flow~\cite{helbing2001traffic,orosz2010traffic}. The subject of traffic dynamics has attracted research interest from many disciplines, including mathematics, physics, and engineering. Since the 1930s, a wide range of models at both the macroscopic and microscopic levels have been proposed to describe traffic behavior~\cite{helbing2001traffic}. Based on these traffic models, many control methods have been introduced and implemented to improve the performance of road transportation systems~\cite{orosz2010traffic}. Currently, most control strategies rely on actuators at fixed locations. For example, variable speed advisory or variable speed limits~\cite{smulders1990control} are commonly implemented through traffic signs on roadside infrastructure, and ramp metering~\cite{papageorgiou2008effects} typically relies on traffic signals located at the freeway entrances. These strategies are essentially \textit{external} regulation methods imposed on traffic flow.

\begin{figure*}[t]
	
	\centering
	\subfigure[]
	{ \label{Fig:TrafficWave1}
		\includegraphics[scale=0.42]{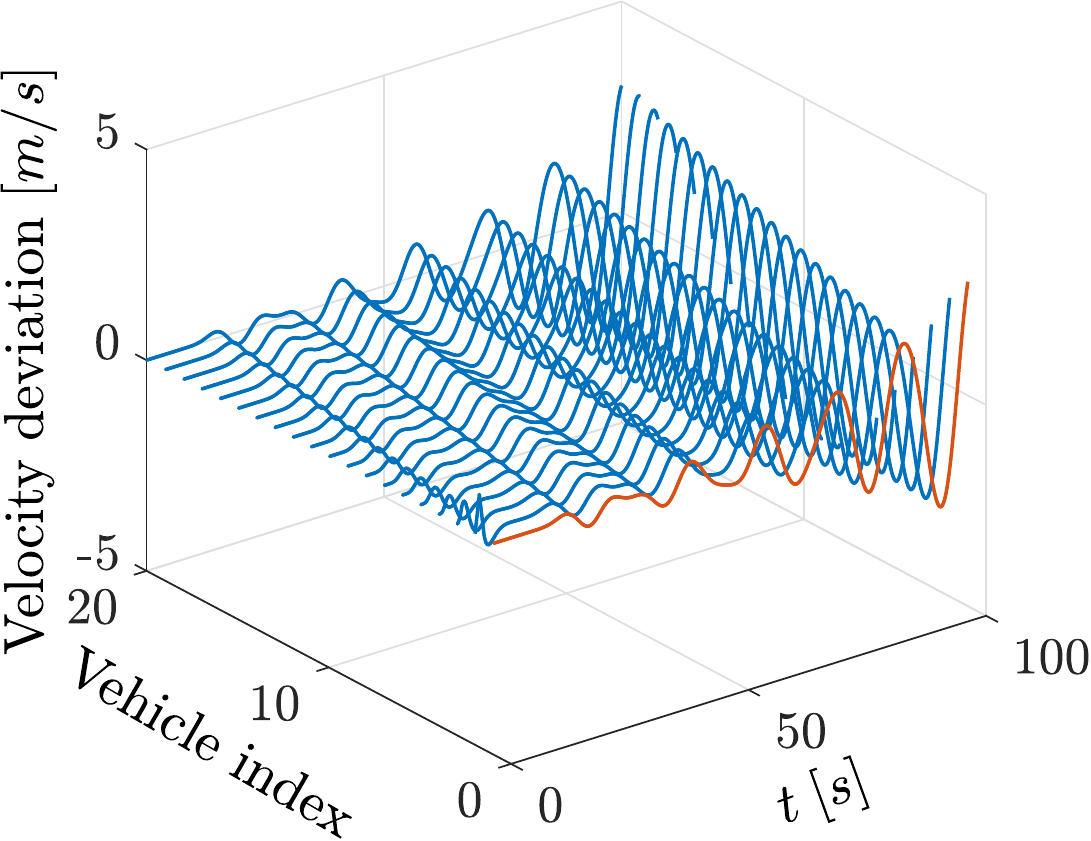}
	}
\hspace{5mm}
	\subfigure[]
	{ \label{Fig:TrafficWave2}
		\includegraphics[scale=0.42]{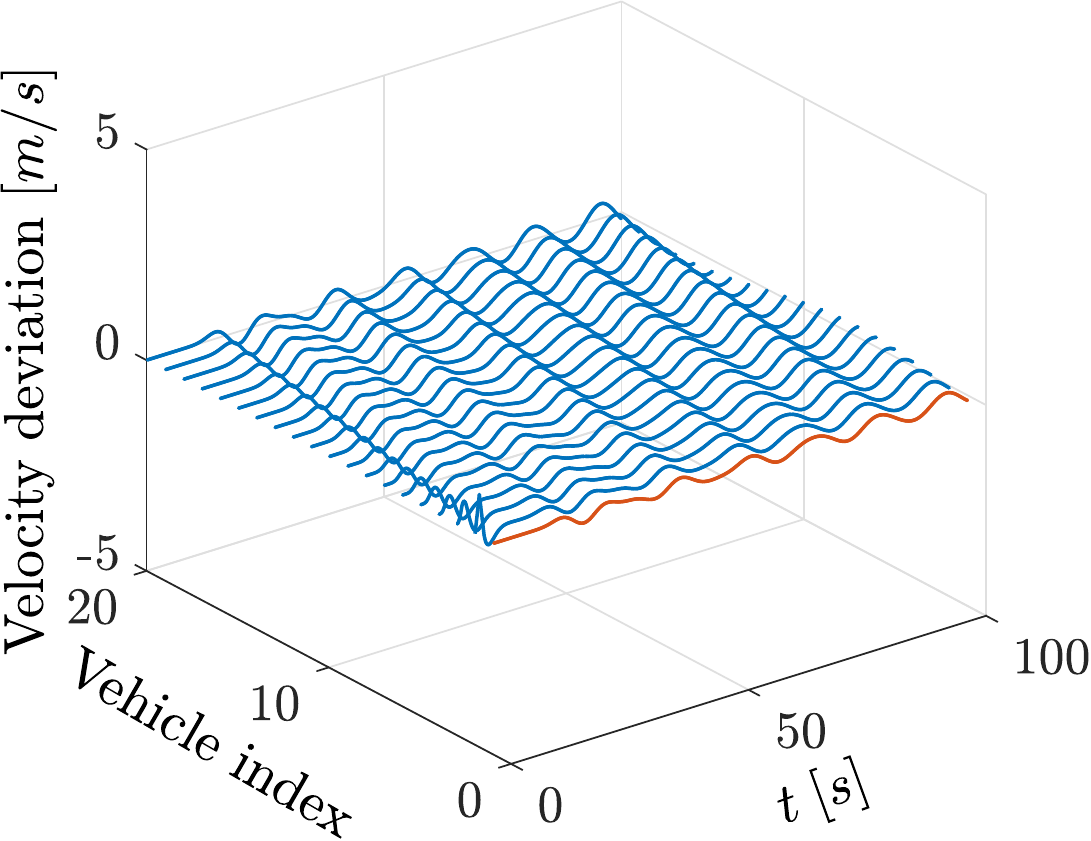}
	}
\hspace{5mm}
	\subfigure[]
	{ \label{Fig:TrafficWave3}
		\includegraphics[scale=0.42]{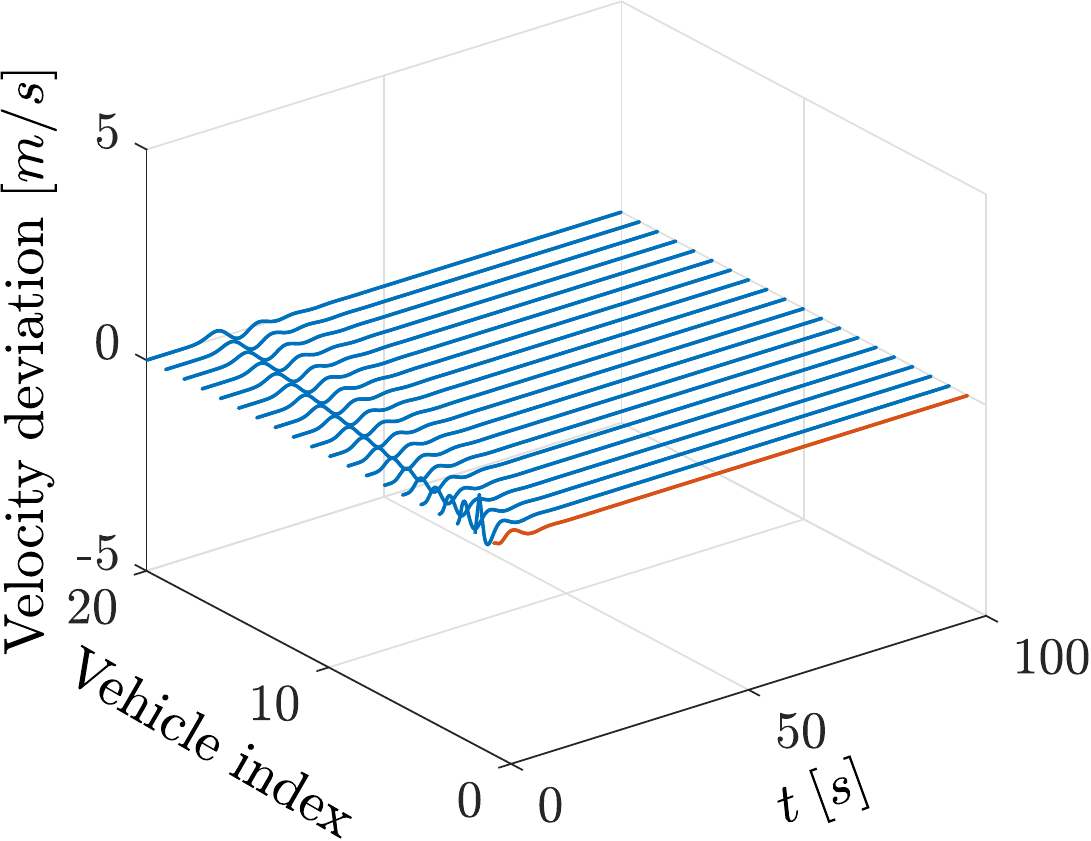}
	}
	\caption{Response profiles to an impulse perturbation in traffic systems on a ring road. Vehicle no.2 has an initial perturbation, and the parameters of human-driven vehicles are chosen to resemble the wave behavior in the real-world experiment~\cite{sugiyama2008traffic}. (a) All the vehicles are human-driven, where the perturbation is amplified and stop-and go waves will appear accordingly. (b) Vehicle no.1 is a CACC-equipped vehicle which adjusts its behavior passively according to its direct preceding vehicle. In this case, the perturbation is not amplified but small traffic waves still persist for a long period. (c) Vehicle no.1 adopts an optimal control strategy considering the global behavior of the entire mixed traffic flow 
to mitigate undesirable perturbations actively. In this case, the perturbation is attenuated, and the traffic flow becomes smooth quickly.}
	\label{Fig:TrafficWave}
\end{figure*}

As a key ingredient of traffic flow, the motion of vehicles plays a fundamental role in road transportation systems. In the past decades, major car-manufacturers and technology companies have invested in developing vehicles with high levels of automation, and some prototypes of {self-driving cars} have been tested in real traffic environments~\cite{fagnant2015preparing}. The emergence of autonomous vehicles (AVs) and vehicular network is expected to revolutionize road transportation~\cite{fagnant2015preparing,contreras2017internet,chen2018driver}. In particular, the advancements of autonomous vehicles offer new opportunities for traffic control, where autonomous vehicles can {receive information from other traffic participants and} act as {mobile} actuators to influence traffic flow \textit{internally}. Most research on the control of traffic flow via autonomous vehicles has focused on platooning of a series of adjacent vehicles or cooperative adaptive cruise control (CACC)~\cite{shladover1991automated,li2017dynamical}. In the context of platoon control, all involved vehicles are assumed to be autonomous and can be controlled to maintain a string stable platoon, such that disturbances along the platoon are dissipated. Significant theoretical and practical advances have been made in designing sophisticated controllers at the platoon level~\cite{naus2010string,milanes2014cooperative,zheng2017distributed}.

While traffic systems with fully autonomous vehicles may be of great interest in the far future, the near future will have to meet a mixed traffic scenario where both autonomous and human-driven vehicles (HDVs) coexist. In fact, early autonomous vehicles need to cooperate in traffic systems where most vehicles are human-driven. This situation is more challenging in terms of theoretical modeling and stability analysis, and many existing studies are based on numerical simulations~\cite{van2006impact,shladover2012impacts,liu2018modeling}. One recent concept is the connected cruise control that considers mixed traffic scenarios where autonomous vehicles can use the information from multiple HDVs ahead to make control decisions~\cite{orosz2016connected,jin2018connected}. More recently, Cui \textit{et al.} first pointed out the potential of a single autonomous vehicle in stabilizing mixed traffic flow~\cite{cui2017stabilizing}, and they also implemented simple control strategies to demonstrate the dissipation of stop-and-go waves via a single autonomous vehicle in real-world experiments~\cite{stern2018dissipation}. The control principle is essentially a slow-in fast-out approach, which is an intuitive method to dampen traffic jams~\cite{nishi2013theory}. This idea was further studied by analyzing head-to-tail string stability~\cite{wu2018stabilizing}. More sophisticated strategies, such as deep reinforcement learning, have also recently been investigated to improve traffic flow in mixed traffic scenarios via numerical simulations~\cite{wu2017flow,wu2017emergent}.

Although the potential of autonomous vehicles has been recognized and demonstrated~\cite{cui2017stabilizing,stern2018dissipation, wu2018stabilizing, wu2017flow, wu2017emergent}, a comprehensive theoretical understanding 
is still lacking. In principle, the behavior of traffic flow emerges from the collective dynamics of many individual human-driven and/or autonomous vehicles~\cite{helbing2001traffic}, where autonomous vehicles can serve as controllable nodes. In this paper, we introduce a control-theoretic framework for analyzing the mixed traffic system with multiple autonomous vehicles and human-driven vehicles. By viewing the autonomous vehicles as controllable nodes, we analyze the controllability and stabilizability of the mixed traffic system. Moreover, we formulate the problem of designing a control strategy for autonomous vehicles to smooth mixed traffic flow as standard $\mathcal{H}_2$ optimal control, and discuss the reachability of desired traffic state. Specifically, the contributions of this paper are:

\begin{enumerate}
\item
We prove {for the first time} that a mixed traffic system consisting of one autonomous vehicle and multiple human-driven vehicles is not completely controllable, but is stabilizable. This result confirms the high potential of autonomous vehicles in mixed traffic systems: the global traffic velocity can be guided to a desired value by controlling autonomous vehicles properly. Our theoretical results validate the empirical experiments in~\cite{stern2018dissipation} that a single autonomous vehicle is able to stabilize traffic flow.  Note that unlike~\cite{cui2017stabilizing}, we prove that there always exists an uncontrollable mode.

\item
We propose an optimal strategy which utilizes a system-level objective for autonomous vehicles: instead of responding to traffic perturbations \textit{passively}, it considers the global behavior of the entire mixed traffic flow to mitigate undesirable perturbations \textit{actively} (see Fig.\ref{Fig:TrafficWave} for illustration). The feedback control of autonomous vehicles is formulated into the standard $\mathcal{H}_2$ optimal synthesis problem~\cite{skogestad2007multivariable}, where the optimal controller can be computed efficiently via convex optimization. Numerical results show a better performance of the proposed controller in smoothing traffic flow and improving fuel economy than existing strategies~\cite{stern2018dissipation}.

\item
Based on reachability analysis, we show an explicit range for the desired traffic state in the mixed traffic system. Moreover, an upper bound of reachable traffic velocity is derived, indicating that a single autonomous vehicle can not only smooth traffic flow but also increase traffic speed. In the numerical experiments, we observed 6\% traffic velocity improvement using only 5\% autonomous vehicles.

\item
Finally, we extend our theoretical framework to the case where multiple autonomous vehicles coexist. We prove that the results on controllability, stabilizability and reachability remain similar. Unlike previous works focusing on one single autonomous vehicle only~\cite{stern2018dissipation,cui2017stabilizing}, we show that autonomous vehicles can cooperate with each other to reduce the time and energy for attenuating perturbations and smoothing traffic flow. As expected, we confirm that controlling multiple autonomous vehicles has a better performance for large-scale mixed traffic systems.
\end{enumerate}

The rest of this paper is organized as follows. Section \ref{Sec:2}¢ò introduces the theoretical modeling of a mixed traffic system with one autonomous vehicle. The controllability and stabilizability result is presented in Section \ref{Sec:3}, and Section \ref{Sec:4} describes the proposed system-level optimal controller and analyzes the reachability of the desired system state. The case of multiple autonomous vehicles is presented in Section \ref{Sec:5}. Numerical simulations are presented in Section \ref{Sec:6}, and we  conclude the paper in Section \ref{Sec:7}.

\section{Theoretical Modelling Framework of Mixed Traffic Systems}
\label{Sec:2}

In this section, we introduce the modeling of a mixed traffic system with one single autonomous vehicle. We consider a single-lane ring road of length $L$ and with $n$ vehicles. As discussed in Cui \textit{et al.}~\cite{cui2017stabilizing}, the ring road setting has several theoretical advantages for modeling a traffic system, including 1) the existence of experimental results that can be used to calibrate model parameters~\cite{sugiyama2008traffic}, 2) perfect control of average traffic density, and 3) correspondence with an infinite straight road with periodic traffic dynamics.

We denote the position of the $i$-th vehicle as $p_i (t)$ along the ring road, and its velocity is denoted as $v_i (t) = \dot{p}_i (t)$. The spacing of vehicle $i$, \textit{i.e.}, the distance between two adjacent vehicles, is defined as $s_i(t) = p_{i-1}(t)-p_i (t)$. Note that we ignore the vehicle length without loss of generality. For simplicity, we assume that there is one autonomous vehicle and the rest are HDVs in this section. The autonomous vehicle is indexed as no.1. The case of multiple autonomous vehicles will be discussed in Section~\ref{Sec:5}.

\subsection{Modeling Human-driven Vehicles}
Human car-following dynamics are typically modeled by nonlinear processes~\cite{helbing2001traffic,orosz2010traffic}
\begin{equation}\label{Eq:HDVNonlinearModel}
	 \dot{v}_i (t)=F(s_i (t),\dot{s}_i (t),v_i (t)),
\end{equation}
meaning that the acceleration of an HDV is a function of its spacing $s_{i}(t)$, the relative velocity between its own and its preceding vehicle $\dot{s}_{i}(t)$, and its velocity $v_{i}(t)$. Denote $s^*,v^*$ as the equilibrium spacing and velocity of each HDV, and then $
(s^{*},v^{*})$ satisfies the following equilibrium equation
\begin{equation}  \label{Eq:EquilibriumEquation}
F(s^{*},0,v^{*})=0,
\end{equation}
which implies a certain relationship between the equilibrium spacing and equilibrium velocity for HDVs. We define the error state of the $i$-th HDV as
\begin{equation*}
\begin{cases}
\tilde{s}_i(t)=s_i(t)-s^*,\\
\tilde{v}_i(t)=v_i(t)-v^*.\\
\end{cases}
\end{equation*}
Applying the first-order Taylor expansion to \eqref{Eq:HDVNonlinearModel} at $(s^{*},v^{*})$ yields a linearized model of each HDV
\begin{equation}\label{Eq:LinearHDVModel}
\begin{cases}
\dot{\tilde{s}}_i(t)=\tilde{v}_{i-1}(t)-\tilde{v}_i(t),\\
\dot{\tilde{v}}_i(t)=\alpha_{1}\tilde{s}_i(t)-\alpha_{2}\tilde{v}_i(t)+\alpha_{3}\tilde{v}_{i-1}(t),\\
\end{cases}
\end{equation}
with $\alpha_1 = \frac{\partial F}{\partial s}, \alpha_2 = \frac{\partial F}{\partial \dot{s}} - \frac{\partial F}{\partial v}, \alpha_3 = \frac{\partial F}{\partial \dot{s}}$  evaluated at the equilibrium state $(s^*,v^*)$. These three coefficients reflect the driver's sensitivity to the error state. Considering the real driver behavior, the acceleration should increase when the spacing increases, the velocity of the ego vehicle drops, or the velocity of the preceding vehicle increases. Hence, we assume that $\alpha _{1}>0$, $\alpha _{2}> \alpha _{3}>0$~\cite{cui2017stabilizing,wilson2011car}.

\begin{figure} [t]
	\centering
	\includegraphics[scale=0.6]{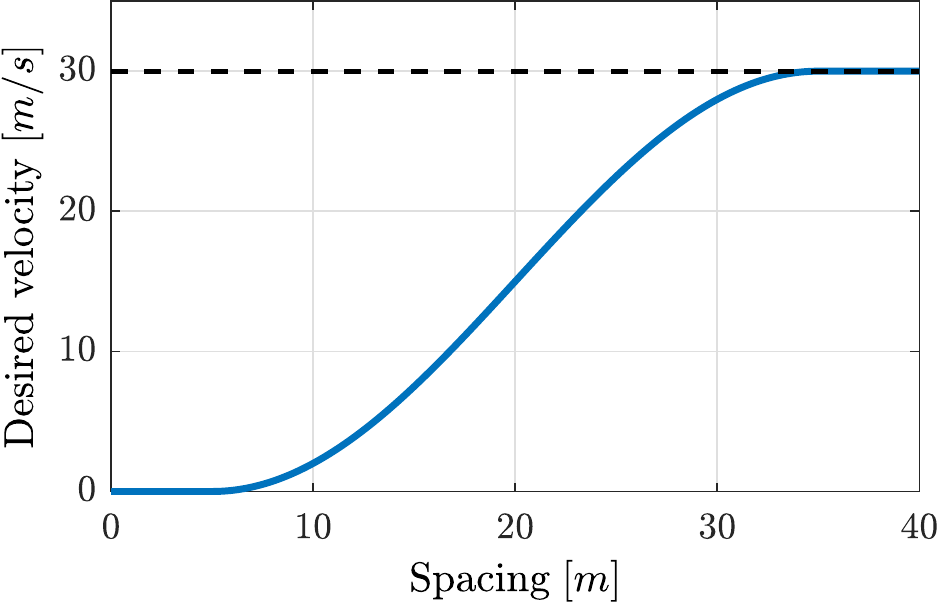}
	\caption{Nonlinear spacing-dependent desired velocity function $V(s)$ in the OVM model. $v_{\text{max}}=30m/s$, $s_{\text{st}}=5m$ and $s_{\text{go}}=35m$. The specific mathematical expression is shown in~\eqref{Eq:OVM_DesiredVelocity} and~\eqref{Eq:OVM_SpacingPolicy}. Note that this figure also illustrates the relationship between equilibrium spacing $s^*$ and equilibrium velocity $v^*$, as shown in \eqref{Eq:OVMequilibrium}.}
	\label{Fig:OVMSpacingPolicy}
\end{figure}

In the following, we choose the optimal velocity model (OVM)~\cite{orosz2010traffic,bando1995dynamical} to derive explicit expressions of~\eqref{Eq:EquilibriumEquation} and~\eqref{Eq:LinearHDVModel}. The OVM model is given by
\begin{equation}   \label{Eq:OVM}
F(s_{i}(t),\dot{s}_{i}(t),v_{i}(t))=\alpha (V(s_{i}(t))-v_{i}(t))+\beta
\dot{s}_{i}(t),
\end{equation}
where $\alpha >0 $ reflects the driver's sensitivity {to} the difference between the current velocity and the spacing-dependent desired velocity $V(s_{i}(t))$, and $\beta > 0$ reflects the driver's sensitivity {to} the difference between the velocities of the ego vehicle and the preceding vehicle. $V(s_{i}(t))$ is usually modeled by a continuous piecewise function
\begin{equation}
\label{Eq:OVM_DesiredVelocity}
V(s)=\begin{cases}
0, &s\le s_{\text{st}},\\
f_v(s), &s_{\text{st}}<s<s_{\text{go}},\\
v_{\max}, &s\ge s_{\text{go}},
\end{cases}
\end{equation}
where the desired velocity $V(s)$ is zero for small spacing $s_{\text{st}}$, and reaches a maximum value $v_{\max }$ for large spacing $s_{\text{go}}$. $f_{v}(s)$ is a monotonically increasing function and defines the desired velocity when the spacing $s$ is between $s_{\text{st}}$ and $s_{\text{go}}$. There are many choices of $f_{v}(s)$, either in a linear or nonlinear form. A typical one is of the following nonlinear form
\begin{equation}
\label{Eq:OVM_SpacingPolicy}
f_{v}(s)=\frac{v_{\max }}{2}\left(1-\cos (\pi
\frac{s-s_{\text{st}}}{s_{\text{go}}-s_{\text{st}}})\right).
\end{equation}
Fig.\ref{Fig:OVMSpacingPolicy} demonstrates a typical example of $V(s)$.

\begin{figure*}[t]
	\centering
	\setlength{\abovecaptionskip}{1em}
	\setlength{\belowcaptionskip}{0.em}
	\includegraphics[scale=0.33]{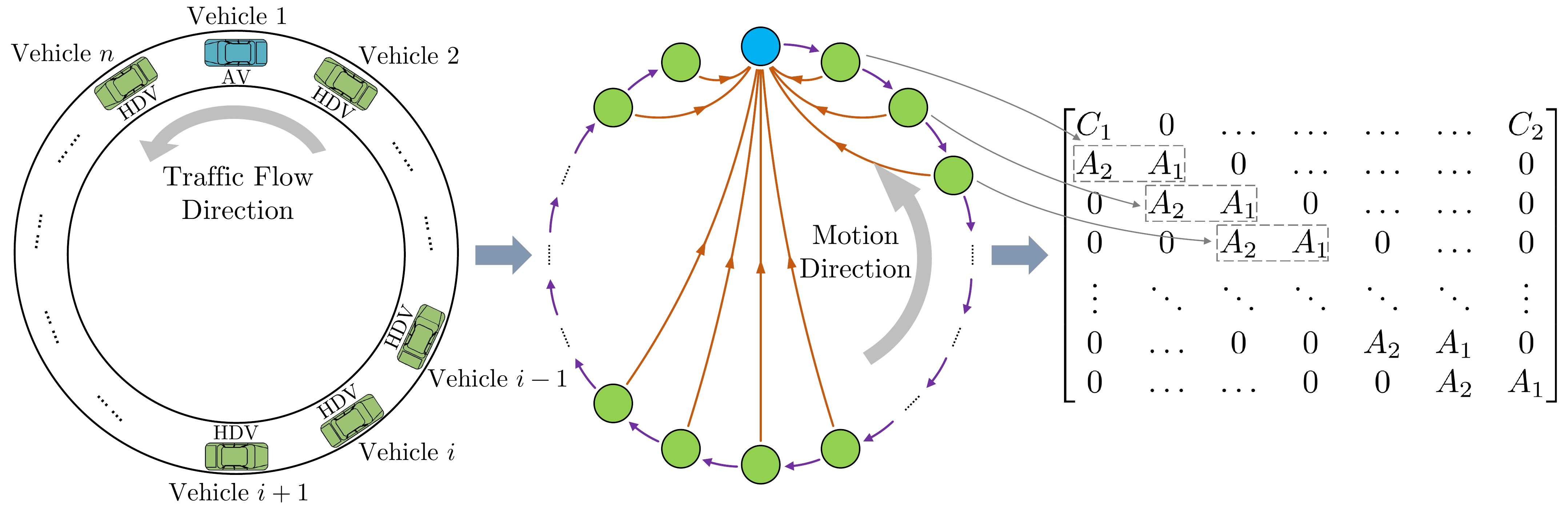}
	\caption{Model establishment schematic. AV: autonomous vehicle; HDV: human-driven vehicle. From left to right,  we have (a) the ring road traffic scenario that includes one autonomous vehicle (blue) and $n-1$ HDVs (green); (b) a simplified network system schematic. Purple arrows indicate the interaction between adjacent vehicles, meaning that each HDV considers the state of its preceding vehicle only. Orange arrows show the information flow of the whole system, assuming that the traffic state is observable to the autonomous vehicle; (c) the system matrix $A$ of the mixed traffic dynamics, as shown in~\eqref{Eq:LinearBlockA}.}
	\label{Fig:SystemModel}
\end{figure*}

For the general OVM model~\eqref{Eq:OVM}, it is easy to obtain the following specific equilibrium state $(s^*,v^*)$ that satisfies~\eqref{Eq:EquilibriumEquation}
\begin{equation}  \label{Eq:OVMequilibrium}
v^* = V(s^*).
\end{equation}
Furthermore, we can calculate the values of the coefficients in linearized model~\eqref{Eq:LinearHDVModel} as follows
\begin{equation} \label{Eq:OVMParameters}
\alpha_1=\alpha \dot{V}(s^*),\ \alpha_2=\alpha+\beta,\ \alpha_3=\beta,
\end{equation}
where $\dot{V}(s^{*})$ denotes the derivative of $V(s)$ with respect to $s$ evaluated at the equilibrium spacing $s^{*}$.

\subsection{Modeling Mixed Traffic Systems}

For the autonomous vehicle, indexed as $i=1$, the acceleration signal is directly used as the control input $u(t)$, and its car-following model is
\begin{equation}\label{Eq:LinearAVModel}
\begin{cases}
\dot{\tilde{s}}_1(t)=\tilde{v}_{n}(t)-\tilde{v}_1(t),\\
\dot{\tilde{v}}_1(t)=u(t),\\
\end{cases}
\end{equation}
where $\tilde{s}_1(t)=s_1(t)-s_c^*,\tilde{v}_1(t)=v_i(t)-v^* $ with $s_c^*$ being a tunable spacing for the autonomous vehicle at velocity $v^*$. Note that $s_c^*$ is a design parameter for the autonomous vehicle, and $v^*$ is the desired traffic velocity, which the autonomous vehicle attempts to steer the traffic flow towards. How to choose a suitable $s_c^*$ is discussed in Section \ref{Sec:Reachability}.

Upon combining the error states of all the vehicles as the mixed traffic
system state,
$$x(t)=\begin{bmatrix} \tilde{s}_1(t),\tilde{v}_1(t),\ldots ,\tilde{s}_n(t),\tilde{v}_n(t) \end{bmatrix}^{\tr}, $$
we arrive at the following canonical linear dynamics from a global system viewpoint
	\begin{equation}
	\dot{x}(t)=Ax(t)+Bu(t), \label{Eq:LinearDynamics}
	\end{equation}
where
	\begin{equation} \label{Eq:LinearBlockA}
	A = \begin{bmatrix} C_1 & 0 & \ldots &\ldots & 0 & C_2 \\
	A_2 & A_1 & 0 & \ldots & \ldots & 0 \\
	0 & A_2 & A_1 & 0 & \ldots & 0\\
	\vdots & \ddots & \ddots & \ddots & \ddots & \vdots\\
	0 & \ldots & 0 & A_2 & A_1 & 0\\
	0 & \ldots & \ldots & 0 & A_2 & A_1\end{bmatrix}, B = \begin{bmatrix} B_1 \\ B_2 \\B_2 \\ \vdots \\ B_2 \end{bmatrix},
	\end{equation}
with each block matrix given by
\begin{gather*} \label{Eq:LinearBlock}
A_1 = \begin{bmatrix} 0 & -1 \\ \alpha_1 & -\alpha_2 \end{bmatrix}, \;
A_2 = \begin{bmatrix} 0 & 1 \\ 0 & \alpha_3 \end{bmatrix}, \;\\
C_1 = \begin{bmatrix} 0 & -1 \\ 0 & 0 \end{bmatrix}, \;
C_2 = \begin{bmatrix} 0 & 1 \\ 0 & 0 \end{bmatrix},  \;
B_1 = \begin{bmatrix} 0  \\ 1 \end{bmatrix},\;
B_2 = \begin{bmatrix} 0 \\ 0 \end{bmatrix}.
\end{gather*}

In~\eqref{Eq:LinearDynamics}, the evolution of each vehicle's state is determined by its own state and the state of its direct preceding vehicle only; see Fig.\ref{Fig:SystemModel} for an illustration. In the following, we provide a theoretical analysis on the potential of the autonomous vehicle on smoothing the mixed traffic flow and design an optimal control input $u(t)$ for the autonomous vehicle.

\section{Controllability and Stabilizability}
\label{Sec:3}

Several real-world field experiments have shown that traffic waves may easily happen and lead to severe congestions in a ring road traffic system with human-driven vehicles~\cite{sugiyama2008traffic,stern2018dissipation}. Theoretical results were also established that the linearized traffic system is stable if the following condition holds~\cite{cui2017stabilizing}
\begin{equation} \label{Eq:StabilityCondition}
\alpha_2^2-\alpha_3^2-2\alpha_1 \geq 0;
\end{equation}
otherwise, the traffic system may become unstable and small perturbations would cause stop-and-go waves. The condition~\eqref{Eq:StabilityCondition} was first derived in~\cite{cui2017stabilizing} using frequency analysis; the interested reader can refer to Appendix~\ref{appendix:stability} for an alternative proof based on eigenvalue analysis.
For the OVM model \eqref{Eq:OVM}, the criterion \eqref{Eq:StabilityCondition} is reduced to
\begin{equation} \label{Eq:StabilityConditionforOVM}
	\alpha+2\beta \geq 2\dot{V}(s^*),
	\end{equation}
which indicates that to guarantee traffic stability, human drivers should have a quicker response to velocity deviations than the sensitivity of the optimal velocity function with respect to the equilibrium spacing; otherwise stop-and-go waves may happen.

Instead of modifying human drivers' behavior, we reveal in this section that a mixed traffic system can be always stabilized by controlling one single autonomous vehicle. Specifically, we discuss two fundamental concepts, controllability and stabilizability, of the mixed traffic system~\eqref{Eq:LinearDynamics}.

\subsection{Controllability Analysis}
\label{Sec:Controllability}

The controllability of a dynamical system captures the ability to guide the system's behavior towards a desired state using appropriate control inputs, and the system is stabilizable if all uncontrollable modes are stable~\cite{skogestad2007multivariable}. We first present three useful lemmas~\cite{kalman1963mathematical, skogestad2007multivariable}.

\begin{lemma}[Controllability] \label{Lemma:Controllability}
	The linear system $(A,B)$ in~\eqref{Eq:LinearDynamics} is controllable if and only if $$\text{rank} \begin{bmatrix}B, AB, \ldots, A^{2n-1}B\end{bmatrix} = 2n.$$
\end{lemma}

Lemma~\ref{Lemma:Controllability} is the well-known Kalman's controllability rank test~\cite{kalman1963mathematical}, which provides a necessary and sufficient mathematical condition for controllability. However, computing the rank requires all the elements of $(A,B)$ to be known, and it might be numerically unreliable to calculate the rank for large-scale systems. To facilitate analysis, we can apply a certain {linear transformation} and represent the linear system under a different basis,  thus simplifying the system dynamics.

In particular, given a nonsingular $T$, we define a new state $\tilde{x} = T^{-1}x$, leading to the following dynamics
 $$
	\dot{\tilde{x}}(t) = T^{-1}AT\tilde{x}(t) + T^{-1}Bu(t).
	$$
Then, we obtain an equivalent linear system $(T^{-1}AT,T^{-1}B)$.

\begin{lemma}[Invariance under linear transformation] \label{Lemma:InvarianceLinear}
	The linear system $(A, B)$ is controllable if and only if $(T^{-1}AT, T^{-1}B)$ is controllable for every nonsingular $T$.
\end{lemma}

If one can diagonalize the system matrix $A$ via $T^{-1}AT$, then the controllability of $(A,B)$ will be easier to derive. This diagonalization may be nontrivial for its original form $(A,B)$. In this case, we can first apply a certain state feedback to simplify the system dynamics. Specifically, consider a control law $ v(t) = u(t) + Kx(t)$, and we arrive at
 $$
	\dot{x}(t) = (A-BK)x(t) + Bv(t).
	$$

\begin{lemma}[Invariance under state feedback] \label{Lemma:InvarianceStateFeedback}
	The linear system $(A, B)$ is controllable if and only if $(A - BK, B)$ is controllable for every $K$ with compatible dimension.
\end{lemma}

Now, we are ready to present the result on the controllability of the mixed traffic system~\eqref{Eq:LinearDynamics}.
\begin{theorem}
	\label{Theorem:Controllability}
	Consider the mixed traffic system in a ring road with one AV and $n-1$  HDVs given by \eqref{Eq:LinearDynamics}. We have
 \begin{enumerate}
   \item The system is not completely controllable.
   \item There exists one uncontrollable mode corresponding to a zero eigenvalue, and this uncontrollable mode is stable.
 \end{enumerate}
\end{theorem}

\begin{IEEEproof}
	Our main idea is to exploit the invariance of controllability under linear transformation and state feedback. Using a sequence of state feedback and linear transformation, we diagonalize the system, leading to an analytical conclusion on the controllability of~\eqref{Eq:LinearDynamics}. Our procedure is as follows.
	$$
	(A,B)\xrightarrow{\text{state feedback}} (\hat{A},B) \xrightarrow{\text{linear transformation}} (\widetilde{A},\widetilde{B})
	$$
	
	First, we transform system $(A,B)$ into $(\hat{A},B)$ by introducing a virtual input $\hat{u}(t)$, defined as
    $$\hat{u}(t)=u(t)-(\alpha_{1}\tilde{s}_{1}(t)-\alpha _{2}\tilde{v}_{1}(t)+\alpha _{3}\tilde{v}_{n}(t)),$$
    which is the difference between the actual control value and the acceleration value when the vehicle is controlled by a human driver. Then, the state
	space model of $(\hat{A},B)$ becomes
	\begin{equation*}
	\dot x (t)=\hat Ax(t)+B\hat u(t),
	\end{equation*}
	where
	\begin{equation} \label{Eq:Ahat}
	\hat A=\begin{bmatrix} A_1 & 0 & \ldots &\ldots & 0 & A_2 \\
	A_2 & A_1 & 0 & \ldots & \ldots & 0 \\
	0 & A_2 & A_1 & 0 & \ldots & 0\\
	\vdots & \ddots & \ddots & \ddots & \ddots & \vdots\\
	0 & \ldots & 0 & A_2 & A_1 & 0\\
	0 & \ldots & \ldots & 0 & A_2 & A_1\end{bmatrix}.
\end{equation}
	
	According to Lemma \ref{Lemma:InvarianceStateFeedback}, the controllability remains the same between system $(\hat{A},B)$ and the original system $(A,B)$. {Note that $\hat{A}$ is block circulant, and it can be diagonalized by the Fourier matrix $F_n$~\cite{marshall2004formations,olson2014circulant}. We refer the interested reader to Appendix~\ref{Appendix:Circulant} for the Fourier matrix $F_n$ as well as precise definitions and properties of block circulant matrices.

Define $\omega=e^{\frac{2\pi j}{n}}$ where $j=\sqrt{-1}$ denotes the imaginary unit, and define $F_n$ given by~\eqref{Eq:FourierMatrix} in Appendix~\ref{Appendix:Circulant}. }Then, we use the transformation matrix $F_{n}^{*}\otimes I_{2}$ to transform $(\hat{A},B)$ into $(\widetilde{A},\widetilde{B})$, where $F_{n}^*$ denotes the conjugate transpose of $F_{n}$. The new system matrix is
	\begin{equation} \label{SEq:TildeADefinition}
{\widetilde A}=(F_{n}^{*}\otimes I_{2})^{-1} \hat A (F_{n}^{*}\otimes I_{2})=\text{diag}(D_{1},D_{2},\ldots
,D_{n}),
\end{equation}
	where $\otimes$ denotes the Kronecker product, and
	\begin{equation}
		\label{Eq:DExpression}
		\begin{aligned}
			D_{i}&=A_{1}+A_{2}\omega^{(n-1)(i-1)}\\
			&=\begin{bmatrix}
		0 & -1+\omega^{(n-1)(i-1)}\\
		\alpha_1 & -\alpha_2+\alpha_3\omega^{(n-1)(i-1)}
		\end{bmatrix}, i=1, 2, \ldots, n,
		\end{aligned}
	\end{equation}
	and $\text{diag}(D_1,D_2,\ldots ,D_n)$ denotes a block-diagonal matrix with $D_1,D_2,\ldots ,D_n$ on its diagonal. Using the fact that $F_n$ is a unitary matrix, the new state variable $\tilde{x}$ after transformation becomes
	\begin{equation}  \label{SEq:TildeXDefinition}
	\tilde{x}=(F_{n}^{*}\otimes I_{2})^{-1}x=(F_{n}\otimes I_{2})x,
	\end{equation}
	and the new control matrix $\widetilde{B}$ is
	\begin{equation*}
	\widetilde{B}=(F_{n}^{*}\otimes I_{2})^{-1}B
	=\frac{1}{\sqrt{n}}\begin{bmatrix}
	B_1\\B_1\\ \vdots \\B_1\\
	\end{bmatrix}.
	\end{equation*}
	Therefore, the dynamics of $\tilde{x}$ are
	\begin{equation} \label{SEq:LinearTransformationS}
	\begin{aligned}
		\dot{\tilde{x}} &= \widetilde{A}\tilde{x}(t)+\widetilde{B}\hat{u}(t)\\
		&=\begin{bmatrix}
	D_{1} & & & & \\
	& D_{2} & & & \\
	& & \ddots & & \\
	& & & D_{n} & \\
	\end{bmatrix}\tilde{x}(t)+\frac{1}{\sqrt{n}}\begin{bmatrix}
	B_1\\B_1\\ \vdots \\B_1
	\end{bmatrix}\hat{u}(t).
	\end{aligned}
	\end{equation}
	Upon denoting $\tilde{x}(t)=\begin{bmatrix}
	\tilde{x}_{11},\tilde{x}_{12},\tilde{x}_{21},\tilde{x}_{22},\ldots ,\tilde{x}_{n1},\tilde{x}_{n2}\end{bmatrix}
	^\tr$, $(\widetilde{A},\widetilde{B})$ is decoupled into $n$ independent subsystems
	\begin{equation*}
	\begin{aligned}
	\frac{d}{dt}{\begin{bmatrix}
		\tilde{x}_{i1} \\
		\tilde{x}_{i2} \\
		\end{bmatrix}
	}&=D_{i} \begin{bmatrix}
	\tilde{x}_{i1} \\
	\tilde{x}_{i2} \\
	\end{bmatrix}
	+ \begin{bmatrix}
	0 \\
	\frac{1}{\sqrt{n}}\\
	\end{bmatrix}
	\hat{u}(t)\\
	&= \begin{bmatrix}
	0 & -1+\omega ^{(n-1)(i-1)} \\
	\alpha _{1} & -\alpha _{2}+\alpha _{3}\omega ^{(n-1)(i-1)}\\
	\end{bmatrix}
	\begin{bmatrix}
	\tilde{x}_{i1}\\
	\tilde{x}_{i2} \\
	\end{bmatrix}
	+\begin{bmatrix}
	0 \\
	\frac{1}{\sqrt{n}}\\
	\end{bmatrix}
	\hat{u}(t).
	\end{aligned}
	\end{equation*}

It is easy to verify that $\dot{\tilde{x}}_{11}=0$, which means that $\tilde{x}_{11}$ is an uncontrollable component, but remains constant during the dynamic evolution. According to  $\widetilde x=(F_{n}^*\otimes I_{2})^{-1}x$, we know
\begin{equation}  \label{SEq:UncontrollableComponent}
\tilde{x}_{11}=\frac{1}{\sqrt{n}}\left((s_{1}(t)-s_{c}^{*})+\sum_{i=2}^{n}{(s_{i}(t)-s^{*})}\right).
\end{equation}

Note that system $(\tilde{A},\tilde{B})$ is equivalent to system $(\hat{A},B)$ due to the linear transformation. Also, system $(\hat{A},B)$ has the same controllability characteristic as the original system $(A,B)$. Therefore, we conclude that the original system $(A,B)$ is not completely controllable and has at least one uncontrollable component which remains constant, as shown in~\eqref{SEq:UncontrollableComponent}.

{Furthermore, the uncontrollable component $\tilde{x}_{11}$ corresponds to a zero eigenvalue. Since this zero eigenvalue only appears in $D_1$ (see the expression of $D_i$ in \eqref{Eq:DExpression}), its algebraic multiplicity in $\hat{A}$ is one. Hence, we conclude that the uncontrollable mode is stable.}
\end{IEEEproof}

We remark that the uncontrollable mode~\eqref{SEq:UncontrollableComponent} has a clear physical interpretation: the sum of each vehicle's spacing should remain constant due to the ring road structure of the mixed traffic system. Theorem~\ref{Theorem:Controllability} differs from the results of~\cite{cui2017stabilizing} in two aspects: 1) we explicitly point out the existence of the uncontrollable mode~\eqref{SEq:UncontrollableComponent}; 2) our proof exploits the block circulant property of the mixed traffic system.

\subsection{Stabilizability Analysis}

After revealing the uncontrollable component~\eqref{SEq:UncontrollableComponent}, we next prove that the mixed traffic system is stabilizable. We need the following PBH test for our stabilizability analysis.
\begin{lemma}[PBH controllability criterion] \label{Lemma:PBHtest}
	The linear system $(A, B)$ is controllable if and only if $\text{rank}(\lambda I - A, B) = 2n$ for every eigenvalue $\lambda$ of $A$. In addition, $(A,B)$ is uncontrollable if and only if there exists $\omega \neq 0$, such that
	$$
	\omega^{\tr} A = \lambda \omega^{\tr}, \omega^{\tr}B = 0,
	$$
	where $\omega $ is a left eigenvector of $A$ corresponding to $\lambda$, and  $\omega$ corresponds to an uncontrollable mode.
\end{lemma}


\begin{theorem}
\label{Theorem:stabilizability}		
Consider the mixed traffic system in a ring road with one AV and $n-1$  HDVs given by \eqref{Eq:LinearDynamics}. We have
\begin{enumerate}
  \item The controllability matrix $Q_c =\begin{bmatrix}B, AB, \ldots, A^{2n-1}B\end{bmatrix} $ satisfies
  \begin{equation}
			\text{rank}(Q_c)= \begin{cases} 2n -1, & \text{if}\;\; \alpha_1 - \alpha_2\alpha_3 + \alpha_3^2 \neq 0,  \\
			n, & \text{if} \;\; \alpha_1 - \alpha_2\alpha_3 + \alpha_3^2 = 0. \end{cases}
		\end{equation}
\item The mixed traffic system~\eqref{Eq:LinearDynamics} is stabilizable.
\end{enumerate}	
\end{theorem}

\begin{IEEEproof}
	As proved in Theorem~\ref{Theorem:Controllability}, systems $(\widetilde{A},\widetilde{B})$ and $(A,B)$ share the same controllability characteristics. Here, we focus on $(\widetilde{A},\widetilde{B})$, and the main idea is to characterize all the uncontrollable modes and prove that they are all stable.

\begin{itemize}
	\item Case 1: $\alpha_1 - \alpha_2\alpha_3 + \alpha_3^2 \neq 0$.
\end{itemize}

	Since $\widetilde{A}$ in~\eqref{SEq:LinearTransformationS} is block diagonal, we have
	\begin{equation}
	\label{Eq:Eigenvalue}
	\text{det}(\lambda I-\widetilde{A})=\prod_{i=1}^{n}{\text{det}(\lambda I-D_{i})}=0,
	\end{equation}
	where $\lambda$ denotes an eigenvalue of $\widetilde{A}$. Substituting~\eqref{Eq:DExpression} into \eqref{Eq:Eigenvalue} leads to the following equation ($i=1,2,\ldots,n$)	
		\begin{equation}
		\label{Eq:EigenvalueEquation}
		\lambda^2+\left(\alpha_2-\alpha_3 \omega^{(n-1)(i-1)}\right)\lambda+\alpha_1\left(1-\omega^{(n-1)(i-1)}\right)
		=0.
		\end{equation}
Note that \eqref{Eq:EigenvalueEquation} is a second-order complex equation and it is non-trivial to directly get its analytical roots. Instead, we use this equation to analyze the properties of the eigenvalues. The following proof is divided into two steps.

	\emph{Step 1:} We prove that $D_i$ and $D_j$ ($i \neq j$) share no common eigenvalues. Assume there exists a $\lambda$ satisfying $\text{det}(\lambda I-D_{i})=0$ and $\text{det}(\lambda I-D_{j})=0$, $i \neq j$, which means
		$$
		\begin{cases}
		\lambda^2+\alpha_2\lambda+\alpha_1=(\alpha_3\lambda+\alpha_1)\omega^{(n-1)(i-1)},\\
		\lambda^2+\alpha_2\lambda+\alpha_1=(\alpha_3\lambda+\alpha_1)\omega^{(n-1)(j-1)}.
		\end{cases}
		$$	
	Since $\omega^{(n-1)(i-1)} \neq \omega^{(n-1)(j-1)}$, we obtain $\alpha_3\lambda+\alpha_1=0$ and $\lambda^2+\alpha_2\lambda+\alpha_1=0$, leading to
		$$
		\alpha_1 - \alpha_2\alpha_3 + \alpha_3^2 = 0,\;\lambda=\alpha_3-\alpha_2,
		$$	
	which contradicts the condition that $\alpha_1 - \alpha_2\alpha_3 + \alpha_3^2 \neq 0$. Therefore, $D_i$ and $D_j (j  \neq i) $ have different eigenvalues.
	
	\emph{Step 2:} We prove that all the system modes corresponding to non-zero eigenvalues are controllable. Denote $\lambda_k \neq 0$ as the eigenvalue of $D_k$ and $\rho$ as its corresponding left eigenvector. According to Lemma~\ref{Lemma:PBHtest}, we need to show $\rho^\tr \widetilde{B} \neq 0 $.
	
	Upon denoting $\rho = \begin{bmatrix}\rho_1^\tr,\rho_2^\tr,\ldots,\rho_n^\tr \end{bmatrix}^\tr$ where $\rho_i=\begin{bmatrix}
	\rho_{i1},\rho_{i2}\end{bmatrix}^\tr \in \mathbb{R}^{2\times1}, i=1,2,\ldots,n $,  the condition $\rho^\tr \widetilde{A} =\lambda_k \rho^\tr $ leads to
		\begin{equation} \label{SEq:Eigenvector}
		\rho_i^\tr D_i = \lambda_k \rho_i^\tr,\;i=1,2,\ldots,n.
		\end{equation}	
	Since $\lambda_k$ is not an eigenvalue of $D_i, \; i  \neq k $, we obtain $\rho_i=0,\; i  \neq k $. Hence, $\rho^\tr \widetilde{B} = \rho_k^\tr B_1 = \rho_{k2}$. Assume $\rho_{k2} = 0$,
	then substituting~\eqref{Eq:DExpression} into~\eqref{SEq:Eigenvector} yields
		\begin{equation}\label{SEq:Rho2Neq0}
		\begin{bmatrix}
		\rho_{k1}&0
		\end{bmatrix}\begin{bmatrix}
		0 & -1+\omega^{(n-1)(k-1)}\\
		\alpha_1 & -\alpha_2+\alpha_3\omega^{(n-1)(k-1)}
		\end{bmatrix}=\lambda_k\begin{bmatrix}
		\rho_{k1}&0
		\end{bmatrix}.
		\end{equation}
	
The only solution to~\eqref{SEq:Rho2Neq0} is $\rho_{k1}=0$, indicating that the left eigenvector $\rho = 0$, which is false. Accordingly, the assumption that $\rho_{k2} = 0 $ does not hold. Therefore, we have $\rho^\tr \widetilde{B} = \rho_{k2} \neq 0 $, meaning that the mode corresponding to $\lambda_k $ is controllable. In other words, the system modes corresponding to nonzero eigenvalues are all controllable. Because the only zero eigenvalue $\lambda=0$ appears in $\text{det}(\lambda I- D_1)=0 $ and  the corresponding mode is uncontrollable, we conclude that if $\alpha_1 - \alpha_2\alpha_3 + \alpha_3^2 \neq 0$, there are $2n-1$ controllable modes in the system $(\widetilde{A}, \widetilde{B})$, meaning that $\text{rank}(Q_c)=2n-1$.
	
	\begin{itemize}
	\item Case 2: $\alpha_1 - \alpha_2\alpha_3 + \alpha_3^2 = 0$.
\end{itemize}

Substituting this condition into \eqref{Eq:EigenvalueEquation} yields
		$$
	\left(\lambda+\alpha_2-\alpha_3\right)\left(\lambda+\alpha_3-\alpha_3\omega^{(n-1)(i-1)}\right)\\
		=0, \; i=1,2,\ldots,n,
		$$	
	which gives the eigenvalues of $D_i$ as follows
		$$
		\lambda_{i1}=\alpha_3-\alpha_2,\;\lambda_{i2}=\alpha_3\left(\omega^{(n-1)(i-1)}-1\right).
		$$
		The rest proof is organized into two steps.
		
	\emph{Step 1:} we prove that there are $n-1$ uncontrollable modes corresponding to $\alpha_3-\alpha_2$.  It is easy to see that $\alpha_3-\alpha_2$ is the common eigenvalue for each block $D_i,\;i=1,2,\ldots,n$, \emph{i.e.}, the algebraic multiplicity of $\alpha_3 - \alpha_2$ is $n$.
	We consider its left eigenvector $\rho = \begin{bmatrix}\rho_1^\tr,\rho_2^\tr,\ldots,\rho_n^\tr \end{bmatrix}^\tr$. Similar to~\eqref{SEq:Eigenvector}, we obtain
		\begin{equation*} \label{SEq:EigenvectorSpecial}
		\rho_i^\tr D_i = (\alpha_3-\alpha_2) \rho_i^\tr,\;i=1,2,\ldots,n.
		\end{equation*}	
	Expanding this equation leads to
		$$
		\begin{bmatrix}
		\rho_{i1}&\rho_{i2}
		\end{bmatrix}\begin{bmatrix}
		0 & -1+\omega^{(n-1)(i-1)}\\
		\alpha_1 & -\alpha_2+\alpha_3\omega^{(n-1)(i-1)}
		\end{bmatrix}=(\alpha_3-\alpha_2)\begin{bmatrix}
		\rho_{i1}&\rho_{i2}
		\end{bmatrix},
		$$	
	from which we have
		$
		\rho_{i1}=-\alpha_3 \rho_{i2},\;i=1,2,\ldots,n.
		$
	Therefore, we can choose $n$ linearly independent left eigenvectors corresponding to $\alpha_3-\alpha_2$ as
		\begin{equation}
		\begin{aligned}
		&\rho^{(1)}=\begin{bmatrix}	-\alpha_3,1,0,0,0,0,\ldots,0,0	\end{bmatrix},\\
		&\rho^{(2)}=\begin{bmatrix}	-\alpha_3,1,\alpha_3,-1,0,0,\ldots,0,0	\end{bmatrix},\\
		&\rho^{(3)}=\begin{bmatrix}	-\alpha_3,1,0,0,\alpha_3,-1,\ldots,0,0	\end{bmatrix},\\
		&\qquad \qquad \qquad \qquad \vdots\\
		&\rho^{(n)}=\begin{bmatrix}	-\alpha_3,1,0,0,0,0\ldots,\alpha_3,-1	\end{bmatrix}.\\
		\end{aligned}
		\end{equation}	
	From these left eigenvectors, it is easy to verify that $\left(\rho^{(1)}\right)^\tr \widetilde{B} \neq 0$ and $\left(\rho^{(i)}\right)^\tr \widetilde{B} = 0,\;i=2,3,\ldots,n $, meaning that for $\alpha_3-\alpha_2$, there are $n-1$ uncontrollable modes.
	
	\emph{Step 2:} we consider the rest of eigenvalues, \emph{i.e.} $$ \lambda_{i2}=\alpha_3\left(\omega^{(n-1)(i-1)}-1\right),\;i=1,2,\ldots,n.$$  The zero eigenvalue $\lambda_{12} = 0$ still corresponds to an uncontrollable mode, as shown in~\eqref{SEq:UncontrollableComponent}. We prove that the modes associated with  $ \lambda_{i2}=\alpha_3\left(\omega^{(n-1)(i-1)}-1\right),\;i=2,3,\ldots,n$ are controllable.
	The proof is similar to the case of $\alpha_1 - \alpha_2\alpha_3 + \alpha_3^2 \neq 0$. For $\lambda_{k2} = \alpha_3\left(\omega^{(n-1)(k-1)}-1\right), \;(k\neq 1)$, denote its left eigenvector as $$\rho = \begin{bmatrix}\rho_1^\tr,\rho_2^\tr,\ldots,\rho_n^\tr \end{bmatrix}^\tr,$$
where $\rho_i=\begin{bmatrix}
	\rho_{i1},\rho_{i2}\end{bmatrix}^\tr \in \mathbb{R}^{2\times1}, i=1,2,\ldots,n $. Then we have $\rho_i=0,\; i  \neq k $, since $\lambda_{k2}$ is not an eigenvalue of other blocks $D_i,\;i\neq k $. For $\rho_k=\begin{bmatrix}\rho_{k1},\rho_{k2}\end{bmatrix}^\tr $, we have $\rho_{k2} \neq 0$, which is similar to the argument in~\eqref{SEq:Rho2Neq0}. Therefore, $\rho^\tr \widetilde{B} \neq 0$, meaning that the mode corresponding to $\lambda_{k2}\; (k \neq 1) $ is controllable.
	
	In summary, the eigenvalue $\lambda = \alpha_3-\alpha_2 $ is associated with $n-1$ uncontrollable modes and one controllable mode. {Since $\lambda = \alpha_3-\alpha_2<0 $, the uncontrollable modes are all stable.} The $n-1$ modes associated with $ \lambda_{i2}=\alpha_3\left(\omega^{(n-1)(i-1)}-1\right),\;i=2,3,\ldots,n$ are controllable, and the zero eigenvalue corresponds to an uncontrollable mode. In total, there are $n$ controllable modes in the system $(\widetilde{A}, \widetilde{B})$, meaning that $\text{rank}(Q_c) = n $. Finally, the system $(A,B)$ is stabilizable since all its uncontrollable mode are stable.
\end{IEEEproof}
{
	Theorem~\ref{Theorem:stabilizability} shows that the mixed traffic system~\eqref{Eq:LinearDynamics} always has one uncontrollable mode corresponding to a zero eigenvalue,
 and the rest of modes are either controllable or stable. This result has no requirement on the car-following behavior of other human-driven vehicles or the scale $n$ of the mixed traffic system. By choosing an appropriate control input, one single autonomous vehicle can always stabilize the global traffic flow at an equilibrium traffic velocity.

\section{Optimal Control and Reachability Analysis}
\label{Sec:4}

We have shown that a mixed traffic system with one single autonomous vehicle is always stabilizable. In this section, we proceed to design an optimal control strategy to reject perturbations in the mixed traffic system using standard control theory~\cite{skogestad2007multivariable}. Moreover, we discuss the reachability of the equilibrium traffic state and show that the autonomous vehicle can increase the traffic equilibrium velocity.

\subsection{Optimal Control Formulation and its Solution}
\label{Sec:OptimalControl}

To reflect traffic perturbations, we assume that there exists a disturbance signal $w_i (t)$ in each vehicle's acceleration signal, \emph{i.e.},
$
\dot{\tilde{v}}_i = \alpha_1 \tilde{s}_i(t) - \alpha_2 \tilde{v}_i(t) + \alpha_3 \tilde{v}_{i-1}(t) + w_i(t).
$
 The linearized dynamics of HDVs in~\eqref{Eq:LinearDynamics} become
$$
\dot{x}_i (t)=A_1 x_i (t)+ A_2 x_{i-1}(t)+H_1w_i(t),
$$
 with $H_1=\begin{bmatrix}0,1\end{bmatrix}^\tr$.

Then, we design an optimal control input $u(t)=-Kx(t)$ to minimize the influence of disturbances $w_i(t)$ on the traffic system, where $K \in \mathbb{R}^{1 \times 2n}$ denotes the feedback gain. Mathematically, this can be formulated into the following optimization problem
 \begin{equation}  \label{SEq:H2problemTraffic}
	\begin{aligned}
	\min_{K} \quad & \|G_{zw}\|^2 \\
	\text{subject to} \quad & u = -Kx,
	\end{aligned}
	\end{equation}
where $G_{zw}$ denotes the transfer function from disturbance signal $w(t)=\begin{bmatrix}w_1(t),\ldots,w_n (t)\end{bmatrix}$ to the performance state
\begin{equation*}
	z(t)=\begin{bmatrix}\gamma_s \tilde s_1 (t),\gamma_v \tilde v_1 (t),\ldots,\gamma_s \tilde s_n (t),\gamma_v \tilde v_n (t),\gamma_u u(t)\end{bmatrix}^\tr,
\end{equation*}
with positive weights $\gamma_s>0,\gamma_v>0,\gamma_u>0$, and $\|\cdot\|$ denotes the $\mathcal{H}_2$ norm of a transfer function that captures the influence of disturbances. Note that the performance state can also be written into
 \begin{equation} \label{Eq:Output}
 	z(t) = \begin{bmatrix} Q^{\frac{1}{2}} \\0 \end{bmatrix}x(t) +  \begin{bmatrix} 0 \\R^{\frac{1}{2}} \end{bmatrix}u(t),
 \end{equation}
where $Q^{\frac{1}{2}} = \text{diag}(\gamma_s,\gamma_v,\ldots,\gamma_s,\gamma_v) $ and $R^{\frac{1}{2}} = \gamma_u $  denote the square roots of state and control performance weights, respectively.
	{
}

The optimization problem~\eqref{SEq:H2problemTraffic} is in the standard form of the $\mathcal{H}_2$ optimal controller synthesis~\cite{skogestad2007multivariable}. Here, we briefly present the steps to obtain a convex formulation for~\eqref{SEq:H2problemTraffic}.
%
\begin{lemma}[$\mathcal{H}_2$ norm of a transfer function\cite{skogestad2007multivariable}]\label{lemma:H2norm}
	Given a stable linear system $\dot{x}(t) = Ax(t)+Hw(t), z(t) = Cx(t)$, the $\mathcal{H}_2$ norm of the transfer function from disturbance $w(t)$ to performance signal $z(t)$ can be computed by
	 $$
		\|G_{zw}\|^2 = \inf_{X \succ 0}\{\text{Trace}\big(CXC^{\tr}\big) \mid AX+XA^{\tr}+HH^{\tr} \preceq 0\},
		$$
	where $\text{Trace}(\cdot)$ denotes the trace of a symmetric matrix.
\end{lemma}

When applying state-feedback $u = -Kx$, the dynamics of the closed-loop traffic system become
\begin{equation}
	\begin{aligned}
	\dot{x}(t) &= (A - BK)x(t) + Hw(t), \\
	z(t) &= \begin{bmatrix} Q^{\frac{1}{2}} \\ -R^{\frac{1}{2}}K \end{bmatrix}x(t).
	\end{aligned}
\end{equation}
Using Lemma~\ref{lemma:H2norm} and a standard change of variables $ Z = KX $, the optimal control problem~\eqref{SEq:H2problemTraffic} can be equivalently reformulated as
  \begin{equation*}
	\begin{aligned}
	\min_{X,Z} \quad & {\text{Trace}}(QX)+{\text{Trace}}\left(RZX^{-1}Z^{\tr}\right) \\
	\text{subject to} \quad & (AX-BZ)+(AX-BZ)^{\tr} + HH^{\tr} \preceq 0, \\
	& X \succ 0.
	\end{aligned}
	\end{equation*}

By introducing $Y \succeq ZX^{-1}Z^{\tr}$ and using the Schur complement, 
a convex reformulation
to~\eqref{SEq:H2problemTraffic} is derived as follows.
 \begin{equation} \label{SEq:LMIHnorm3}
	\begin{aligned}
	\min_{X,Y,Z} \quad & \text{Trace}(QX)+{\text{Trace}}(RY)  \\
	\text{subject to} \quad & (AX-BZ)+(AX-BZ)^{\tr} + HH^{\tr} \preceq 0, \\
	& \begin{bmatrix} Y & Z \\ Z^{\tr} & X \end{bmatrix} \succeq 0, X \succ 0.
	\end{aligned}
	\end{equation}
Problem~\eqref{SEq:LMIHnorm3} is convex and ready to be solved using general conic solvers, \emph{e.g.}, Mosek~\cite{mosek2010mosek}, and the optimal controller is recovered as $K = ZX^{-1}$.

\begin{remark}[Active response to traffic perturbations]
 Some traditional strategies for autonomous vehicles, \eg, CACC \cite{naus2010string,milanes2014cooperative}, mainly focus on the performance of the autonomous vehicles themselves, {corresponding to a local-level consideration,} and they typically respond to external perturbations in a passive way. Although they can improve traffic stability in mixed traffic flow~\cite{van2006impact,shladover2012impacts}, there always exists a certain requirement on the penetration rate of autonomous vehicles. Instead, our formulation directly considers the global traffic behavior, \ie, the state and behavior of all the involved vehicles, and aims at minimizing the influence of disturbances on the entire traffic system via controlling autonomous vehicles. In this way, the resulting {system-level} strategy enables autonomous vehicles to respond to traffic perturbations actively; see Fig.~\ref{Fig:TrafficWave} for an illustration. 
\end{remark}

\begin{remark}[Parameter selection of controllers] \label{Remark:Parameter}
{There already exist several control strategies for autonomous vehicles to stabilize traffic flow, \emph{e.g.}, FollowerStopper and PI with Saturation in~\cite{stern2018dissipation}. We note that there are many parameters that need to be pre-determined in these strategies. Different choices of parameter values may lead to different performance, which are not completely predictable.} Instead, only three parameters need to be designed in~\eqref{SEq:H2problemTraffic}, \emph{i.e.}, the weight coefficients in the cost function, $\gamma_s$, $\gamma_v$, and $\gamma_u$. Moreover, we can adjust their values to achieve different and predictable results. For example, setting a larger value to $\gamma_s$ and $\gamma_v$ typically allows to stabilize the traffic in a shorter time, and setting a larger value to $\gamma_u$ normally helps to keep a lower control energy for the autonomous vehicle.
\end{remark}

\begin{remark} [Feasibility of the proposed strategy]
Problem~\eqref{SEq:LMIHnorm3} can be formulated into a standard semidefinite program (SDP), for which there exist efficient algorithms to get a solution of arbitral accuracy in polynomial time~\cite{ben2001lectures}. In particular, the well-established interior-point method has a computational complexity of $\mathcal{O}(\max(m,n) m n^{2.5})$, where $n$ is the dimension of the positive semidefinite constraint and $m$ is the number of equality constraints in a standard SDP~\cite{ben2001lectures}. In this paper, we use the conic solver Mosek~\cite{mosek2010mosek} to solve the resulting SDP from~\eqref{SEq:LMIHnorm3}. In addition, computing the optimal control strategy from~\eqref{SEq:LMIHnorm3} requires explicit car-following models of HDVs; see~\eqref{Eq:HDVNonlinearModel}. Lyapunov-type methods~\cite{gomez2007line} and data-driven methods~\cite{lin2018moha} have been proposed for practical estimation of car-following behaviors based on vehicle trajectories. It would be interesting for future work to incorporate model identification in the design of robust optimal control strategies.
\end{remark}

\subsection{Reachability and Maximum Traffic Velocity}
\label{Sec:Reachability}

We have shown that the mixed traffic system is always stabilizable and upon choosing the weight coefficients $\gamma_s$, $\gamma_v$, and $\gamma_u$, we can solve~\eqref{SEq:LMIHnorm3} to obtain an optimal control strategy $u(t) = - K x(t)$ to reject the influence of disturbances. Considering the definition of state $x(t)$ and denoting $K = \begin{bmatrix} k_{1,1},k_{1,2},k_{2,1},k_{2,2,}, \ldots, k_{n,1},k_{n,2} \end{bmatrix}\in \mathbb{R}^{1\times 2n}$, the optimal control strategy is implemented as
\begin{equation}  \label{SEq:OptimalControl}
\begin{aligned}
u(t) = &- \left(k_{1,1}\left(s_{1}(t)-s^{*}_c\right) +  k_{1,2}\left(v_{1}(t)-v^{*}\right) \right) \\
&- \sum_{i=2}^n\left(k_{i,1}(s_{i}(t)-s^{*}) + k_{i,2}(v_{i}(t)-v^{*})\right).
\end{aligned}
\end{equation}
Recall that $(s^*,v^*)$ is the traffic equilibrium state of HDVs satisfying~\eqref{Eq:EquilibriumEquation}, and $s^*_c>0$ is the desired spacing for the autonomous vehicle that is free to choose. We observe that an improper choice of $s_c^*$ may cause the mixed traffic system to fail to reach the desired velocity $v^*$ due to the uncontrollable mode~\eqref{SEq:UncontrollableComponent}.
In fact, we can further predict the system final state, which sheds some insight on the reachability of the equilibrium traffic state.

The result of reachability is stated as follows.

\begin{theorem}
\label{Theorem:reachability}	
Consider the mixed traffic system in a ring road with one AV and $n-1$  HDVs given by \eqref{Eq:LinearDynamics}. Suppose a stabilizing feedback gain is found by \eqref{SEq:H2problemTraffic} and the matrix~\eqref{eq:Mnonsingular} is non-singular. Then the traffic system is stabilized at $v^*$ if and only if the desired spacing $s_c^*$ of the AV satisfies
	\begin{equation} \label{Eq:ReachabilityCondition}
	s_c^*=L-(n-1) s^*,
	\end{equation}
	with $s^*$ given by the equilibrium equation \eqref{Eq:EquilibriumEquation}.
\end{theorem}

\begin{IEEEproof}
	With a stabilizing controller for the autonomous vehicle, the mixed traffic system~\eqref{Eq:LinearDynamics} is stable and the state $x(t)$ will approach an equilibrium point, where $\dot x (t) =0$. We analyze the dynamics of each vehicle in~\eqref{Eq:LinearHDVModel} and~\eqref{Eq:LinearAVModel} separately, leading to
	\begin{equation*}
	\tilde{s}_{1}(t_{\text{f}})=s_{\text{e}}, \; \tilde{v}_{i}(t_{\text{f}})=v_{\text{e}}, \;
	\tilde{s}_{i}(t_{\text{f}})=\frac{\alpha _{2}-\alpha _{3}}{\alpha _{1}}v_{e}, \; i=2,
	3,\ldots, n,
	\end{equation*}
	where $s_{\text{e}}, v_{\text{e}}$ are constant values, and $t_{\text{f}}$ is the time when the system reaches its equilibrium point. Considering the desired state in the controller $x_{\text{des}}=\begin{bmatrix}
	s_{c}^{*} ,v^{*},s^{*},v^{*},\ldots ,s^{*},v^{*}
	\end{bmatrix}^\tr $, the final state of the
	system~\eqref{Eq:LinearDynamics} must be in the following form
	\begin{equation}
	x_{\text{f}}=\begin{bmatrix}
	s_\text{f,AV},v_f,s_\text{f,HDV},v_f,\ldots ,s_\text{f,HDV},v_f
	\end{bmatrix}^\tr,
	\end{equation}
	where $s_\text{f,AV}= s_{c}^{*}+s_{e}$, $s_\text{f,HDV}=s^{*}+ \displaystyle \frac{\alpha _{2}-\alpha _{3}}{\alpha_{1}} v_{e}$, $v_f=v^{*}+v_{e}$.
	
	We next calculate the exact value of $s_{e}$ and $v_{e}$. In the final state, all the vehicles have zero acceleration, indicating that the control input $u(t)$ must be zero, \textit{i.e.}, $u(t)=-Kx(t)=0$. Besides, according to the
	controllability analysis, we know that there exists an uncontrollable mode $(s_{1}(t)-s_{c}^*)+\sum_{i=2}^{n}{(s_{i}(t)-s^{*})}$ remaining constant; see~\eqref{SEq:UncontrollableComponent}. Combining these two conditions leads to the following linear equations
	\begin{equation} \label{Eq:Reachability}
	\begin{cases}
	\left( \frac{\alpha_2-\alpha_3}{\alpha_1} \displaystyle \Sigma_{i=2}^n k_{i,1} + \displaystyle \Sigma_{i=1}^n k_{i,2}  \right)v_e + k_{1,1} s_e = 0,    \\
	(n-1)\frac{\alpha_2-\alpha_3}{\alpha_1}  v_e  + s_e =L -(n-1)s^*- s_c^*,
	\end{cases}
	\end{equation}
	and its solution offers the exact value of $s_{e}$ and $v_{e}$, \emph{i.e.},
    $$
        \begin{bmatrix} s_e \\ v_e \end{bmatrix} = M^{-1} \begin{bmatrix} 0 \\L -(n-1)s^*- s_c^* \end{bmatrix},
    $$
    where $M$ is a non-singular matrix as
    \begin{equation} \label{eq:Mnonsingular}
       M =  \begin{bmatrix} \frac{\alpha_2-\alpha_3}{\alpha_1} \displaystyle \Sigma_{i=2}^n k_{i,1} + \displaystyle \Sigma_{i=1}^n k_{i,2}  &  k_{1,1} \\
       (n-1)\frac{\alpha_2-\alpha_3}{\alpha_1}  & 1   \end{bmatrix}.
    \end{equation}

	To reach the desired equilibrium state $(s^*,v^* )$, we should have $s_e=0$ and $v_e=0$, which is equivalent to
$
    L -(n-1)s^*- s_c^*=0.
$
 This completes our proof.
\end{IEEEproof}

Since the mixed traffic system in~\eqref{Eq:LinearDynamics} is stabilizable, it can be guided to reach any equilibrium state with traffic speed $v^*$  via controlling the autonomous vehicle properly. In practice, however, the spacing of the autonomous vehicle cannot be negative, \textit{i.e.}, $s_c^*>0$, which is equivalent to
	\begin{equation}\label{Eq:MaximumSpacing}
	s^*_{\max} < \frac{L}{n-1}.
	\end{equation}
Recall that $(s^*,v^*)$ should satisfy the HDV equilibrium equation $F(s^*,0,v^*)$, and $v^*$ usually increases as $s^*$ grows up according to the real driving behavior, as illustrated in Fig.\ref{Fig:OVMSpacingPolicy} for the OVM model. Therefore, the requirement of the equilibrium spacing in \eqref{Eq:MaximumSpacing} sets up a maximum equilibrium traffic velocity $v^*_{\max}$. This leads to the following result.

 \begin{corollary}[Range of reachable traffic velocity] \label{Corollary:Velocity}
 	There exists a reachable range for the traffic velocity in the mixed traffic system~\eqref{Eq:LinearDynamics}:
 	\begin{equation}
 		0\leqslant v^* < v^*_{\max},
 	\end{equation}
 	 where $v^*_{\max}$ satisfies
 	\begin{equation*}
 	F\left(\frac{L}{n-1},0,v^*_{\max}\right)=0.
 	\end{equation*}
 \end{corollary}

 This result reveals the upper bound of reachable traffic velocity, indicating that mixed traffic flow with one autonomous vehicle can be steered exactly towards a velocity within the range of $\left(0, v^*_{\max}\right)$. Note that $v^*_{\max}$ is higher than the equilibrium traffic speed with HDVs only, where the equilibrium spacing is $s^*= L/n$. A physical interpretation is that the autonomous vehicle can follow its preceding vehicle at a shorter distance and leave more space for its following HDVs, which in turn triggers the HDVs to travel at a higher speed in the equilibrium; see Fig.\ref{SFig:IncreasingSpeed} for illustration.

\begin{remark} 
	In the case of vehicle platoons where all involved vehicles have autonomous capabilities, the vehicles can be controlled to reach the same desired velocity with separate desired spacings~\cite{naus2010string,milanes2014cooperative,zheng2017distributed}. In a mixed traffic system, although only the autonomous vehicles are under direct control, all the other HDVs can be influenced indirectly. 
One distinction is that the desired state $(s^*,v^* )$ for HDVs should satisfy their corresponding car-following behaviors, while the desired state $(s_c^*,v^* )$ for the autonomous vehicle can be designed separately.
\end{remark}

\begin{figure}[t]
	\centering
	\setlength{\abovecaptionskip}{-0.2em}
	\setlength{\belowcaptionskip}{0em}
	\newcommand{\fighspace}{\hspace{10mm}}
	\subfigure[]
	{ \label{Fig:IncreasingSpeeda}
		\includegraphics[scale=0.4]{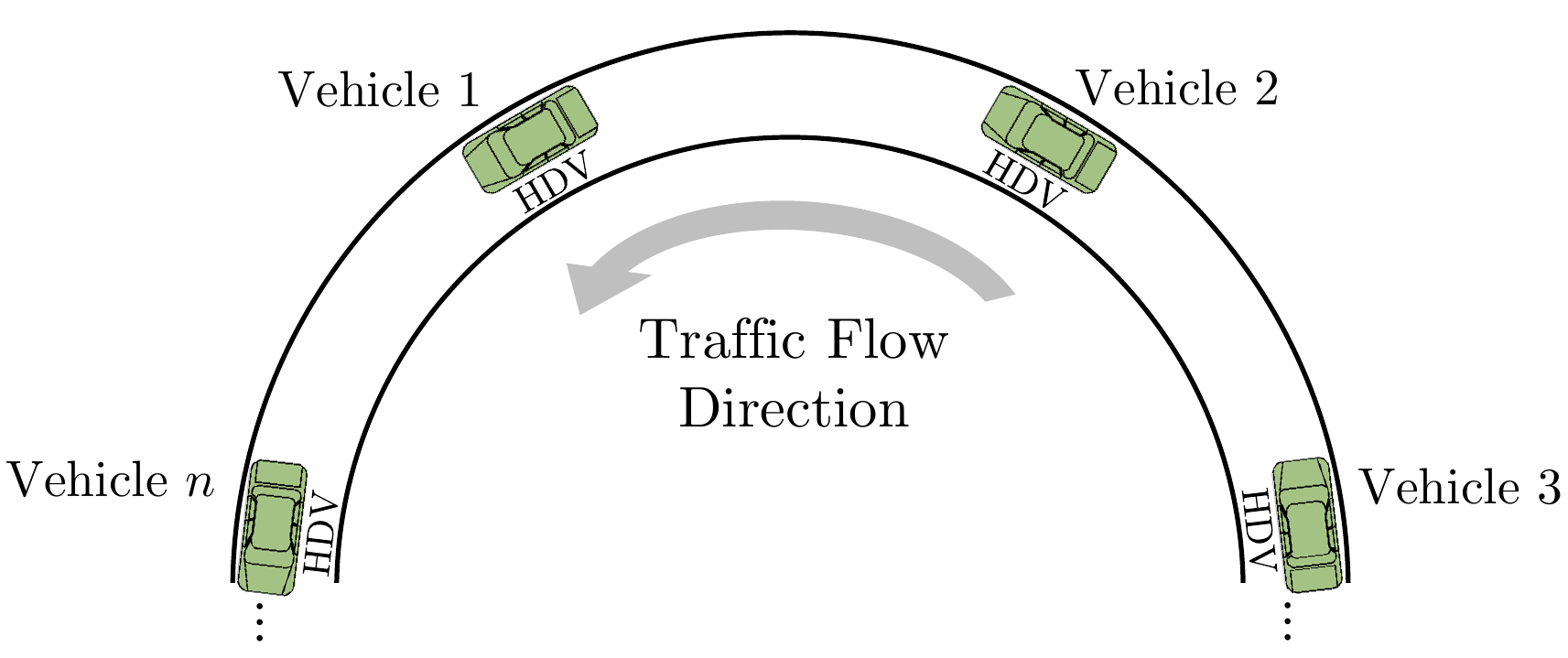}
	} 
	\subfigure[]
	{ \label{Fig:IncreasingSpeedb}
		\includegraphics[scale=0.4]{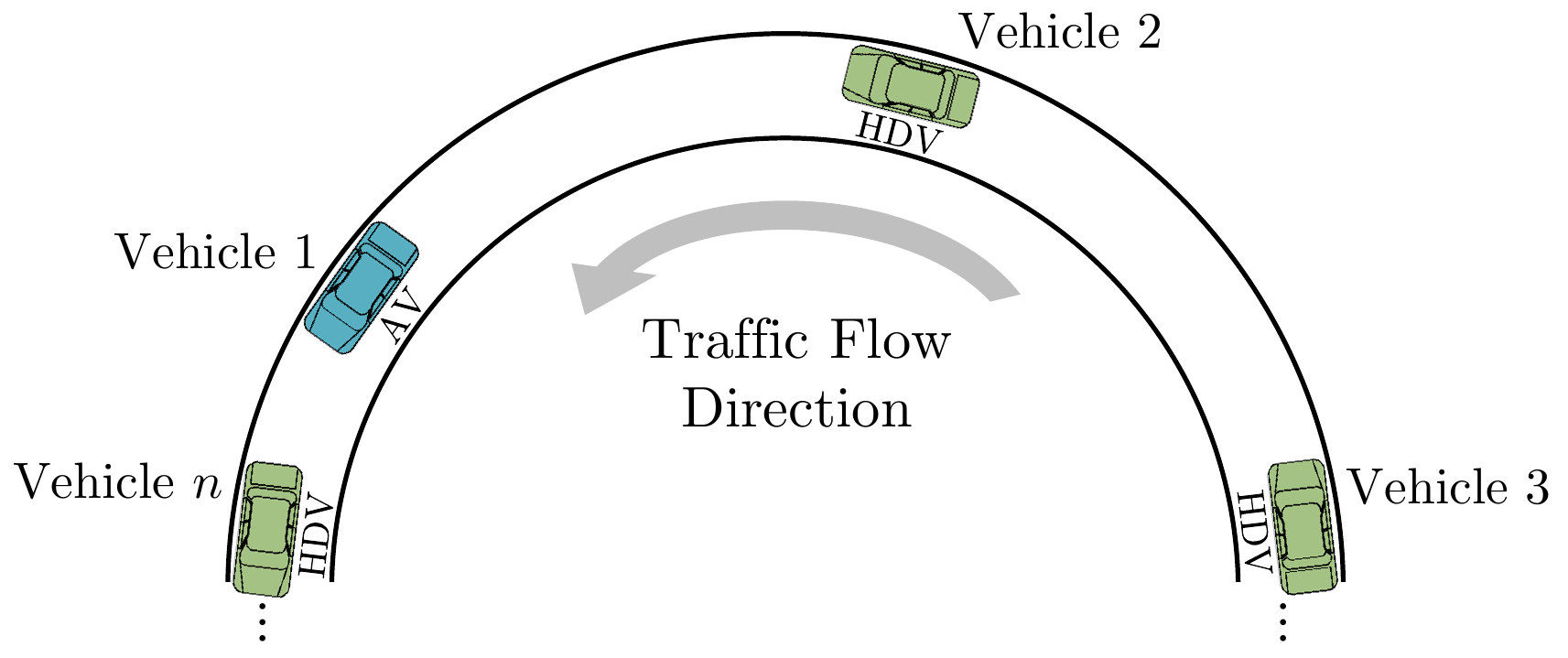}
	}
	\caption{Illustration of the scenario where the autonomous vehicle increases the traffic speed. (a) When all vehicles are human-driven, the spacing between two vehicles is equal for homogeneous car-following dynamics. (b) In the case of mixed traffic systems, the autonomous vehicle can be controlled to follow its preceding vehicle in a shorter distance, and the other HDVs have a larger spacing at the equilibrium state. According to $F(s^*,0,v^*)=0$, the equilibrium velocity $v^*$ increases as $s^*$ grows up. Hence, the entire traffic flow speed can be increased via controlling the autonomous vehicle.}
	\label{SFig:IncreasingSpeed}
\end{figure}

\section{Traffic Systems with Multiple Autonomous Vehicles}
\label{Sec:5}

The mixed traffic system with a single autonomous vehicle is stabilizable, which is independent of the number of vehicles $n$. Also, an optimal control strategy $u(t)$ can be obtained by solving an optimization problem. However, it might be not practical to control a mixed traffic system consisting of many HDVs and a single autonomous vehicle. In this section, we extend our previous analysis to a mixed traffic system with multiple autonomous vehicles.

\subsection{Theoretical Framework}

Assume that there are $n$ vehicles in the traffic flow with $k$ autonomous vehicles $(k<n)$. The indices of the autonomous vehicles are $i_{1}
$, $i_{2}$, \ldots , $i_{k}$, for which we define a set $S_{\text{AV}}=\{ i_{1},i_{2},\ldots ,i_{k}\} $. The error state of an HDV indexed as $i$ ($i\notin S_{\text{AV}}$) is still defined as $\tilde{s}_i(t)=s_i(t)-s^*, \tilde{v}_i(t)=v_i(t)-v^*$, where ($s^*,v^*$) satisfies \eqref{Eq:EquilibriumEquation} and its linearized model remains the same as \eqref{Eq:LinearHDVModel}.

For an autonomous vehicle indexed as $i_r$ ($r=1,2\ldots,k$), the acceleration signal is directly used as its control input $u_{i_r}(t)$, and its car-following model is
\begin{equation}\label{Eq:LinearAVModelMultiple}
\begin{cases}
\dot{\tilde{s}}_{i_r}(t)=\tilde{v}_{{i_r}-1}(t)-\tilde{v}_{i_r}(t),\\
\dot{\tilde{v}}_{i_r}(t)=u_{i_r}(t),\\
\end{cases}
\end{equation}
where $\tilde{s}_{i_r}(t)=s_{i_r}(t)-s_{{i_r},c}^*,\tilde{v}_{i_r}(t)=v_{i_r}(t)-v^* $ with $s_{{i_r},c}^*$ being a tunable desired spacing for the autonomous vehicle at velocity $v^*$.

To derive the global dynamics, we lump the error states of all the vehicles as the mixed traffic
system state,
$$x(t)=\begin{bmatrix} \tilde{s}_1(t),\tilde{v}_1(t),\ldots ,\tilde{s}_n(t),\tilde{v}_n(t) \end{bmatrix}^{\tr}, $$
and lump all the control inputs as
$$
u(t)=\begin{bmatrix} u_{i_{1}}(t),u_{i_{2}}(t),\ldots ,u_{i_{k}}(t) \end{bmatrix}^{\tr}.
$$
Then the state-space model of the entire mixed traffic system can be given by
\begin{equation}  \label{SEq:MultipleAVModel}
\dot x (t)=A_kx(t)+B_ku(t),
\end{equation}
where
\begin{equation*}
	A_k=\begin{bmatrix} A_{11} & 0 & \ldots &\ldots & 0 & A_{12} \\
	A_{22} & A_{21} & 0 & \ldots & \ldots & 0 \\
	0 & A_{32} & A_{31} & 0 & \ldots & 0\\
	\vdots & \ddots & \ddots & \ddots & \ddots & \vdots\\
	0 & \ldots & 0 & A_{(n-1)2} & A_{(n-1)1} & 0\\
	0 & \ldots & \ldots & 0 & A_{n2} & A_{n1}\end{bmatrix},
\end{equation*}
\begin{equation*}
	B_k=\begin{bmatrix} P_1,P_2,\ldots,P_k \end{bmatrix}.
\end{equation*}
In the system matrix $A_{k}$, we have
\begin{equation*}
\begin{cases}
A_{r2}=C_{2} ,A_{r1}=C_{1},& \text{if}\;\;  r\in S_{\text{AV}},\\
A_{r2}=A_{2},A_{r1}=A_{1},&  \text{if}\;\;  r\notin S_{\text{AV}},\\
\end{cases}
\end{equation*}
and the other blocks are zero, where $A_{1}$, $A_{2}$, $C_{1}$ and $C_{2}$ are the same as that in~\eqref{Eq:LinearBlockA}. In $B_{k}$, each column $P_{r}$ is a $2n\times 1$ vector, in which only the {$(2i_{r})$-th} entry is one and the others are zero.
\subsection{Controllability and Stabilizability}

When there exist multiple autonomous vehicles, we observe that the controllability and stabilizability of the mixed traffic system remain unchanged compared to the case discussed in Section \ref{Sec:3}, where there is only one autonomous vehicle. This result is summarized as follows.

\begin{theorem}
    Consider a mixed traffic system with multiple autonomous vehicles in~\eqref{SEq:MultipleAVModel}. We have
    \begin{enumerate}
      \item The mixed traffic system~\eqref{SEq:MultipleAVModel} is not completely controllable, and there still exists an uncontrollable mode.
      \item The mixed traffic system~\eqref{SEq:MultipleAVModel} is stabilizable.
    \end{enumerate}

\end{theorem}

\begin{IEEEproof}
To analyze the system controllability, we first define a virtual control input as
$$
\hat{u}_k(t)=u(t)-\begin{bmatrix}\bar u_{i_1}, \bar u_{i_2}\ldots ,\bar u_{i_k}
\end{bmatrix} ^{\tr},
$$
with $\bar u_{i_r}=\alpha
_{1}\tilde{s}_{i_r}(t)-\alpha _{2}\tilde{v}_{i_r}(t)+\alpha
_{3}\tilde{v}_{i_r-1}(t)$, $i_{r}\in S_{\text{AV}}$. Then~\eqref{SEq:MultipleAVModel} becomes
\begin{equation*}
\dot x (t)=\hat Ax(t)+B_k\hat u_k(t),
\end{equation*}
with $\hat A$ defined by \eqref{Eq:Ahat}.
Using $F_{n}^{*}\otimes I_{2}$ as the transformation matrix, $(\hat{A},B_{k})$
can be transformed into $(\widetilde{A},\widetilde{B}_{k})$, given by
\begin{equation} \label{SEq:LinearTransformationM}
\dot{\tilde{x}} = \widetilde{A}\tilde{x}(t)+{\widetilde{B}_k}\hat{u}_k(t)
\end{equation}
with $\widetilde{A}=\text{diag}(D_1,D_2,\ldots,D_n)$ being the same as~\eqref{SEq:TildeADefinition} and $\widetilde{B}_{k}$ defined as
\begin{equation*}
\widetilde{B}_{k}=(F_{n}^{*}\otimes I_{2})^{-1}B_k
=\begin{bmatrix} \widetilde P_1,\widetilde P_2,\ldots,\widetilde P_k \end{bmatrix},
\end{equation*}
where $\widetilde{P}_{r}=\frac{1}{\sqrt{n}}\begin{bmatrix}
0,1,0,\bar{\omega }^{i_{r}-1},\ldots
,0,\bar{\omega }^{(n-1)(i_{r}-1)}
\end{bmatrix}^{\tr}, r = 1, \ldots, k$. After the transformation, the new state variable $\tilde{x}$ is the same as~\eqref{SEq:TildeXDefinition}.

Note that the dynamics~\eqref{SEq:LinearTransformationM} shall be reduced to the case with a single autonomous vehicle~\eqref{SEq:LinearTransformationS} when $k = 1$ and $i_1 = 1$. Upon denoting $\tilde{x}(t)=\begin{bmatrix} \tilde{x}_{11},\tilde{x}_{12},\tilde{x}_{21},\tilde{x}_{22},\ldots ,\tilde{x}_{n1},\tilde{x}_{n2} \end{bmatrix}^{\tr}$, $(\widetilde{A},\widetilde{B}_{k})$ can be decoupled into $n$ independent
subsystems ($q=1,2,\ldots ,n$)
\begin{equation*}
\frac{d}{dt}{\begin{bmatrix}
	\tilde{x}_{q1} \\
	\tilde{x}_{q2} \\
	\end{bmatrix}
}=D_{i} \begin{bmatrix}
\tilde{x}_{q1} \\
\tilde{x}_{q2} \\
\end{bmatrix}
+F_q\hat{u}(t) ,
\end{equation*}
where
$$
F_q=\frac{1}{\sqrt{n}} \begin{bmatrix}
0 & 0 & \cdots & 0 \\
\bar{\omega }^{(q-1)(i_{1}-1)} & \bar{\omega }^{(q-1)(i_{2}-1)} & \cdots &
\bar{\omega }^{(q-1)(i_{k}-1)} \\
\end{bmatrix}.
$$
It is not difficult to see that $\dot {\tilde {{x}}}_{11}=0$, which means that $\tilde x_{11}$ is an uncontrollable mode. Thus, the mixed traffic system~\eqref{SEq:MultipleAVModel} is not completely controllable. Note that $\tilde x_{11}$ corresponds to the zero eigenvalue and remains constant. Similar to the analysis in Section \ref{Sec:Controllability}, the algebraic multiplicity of the zero eigenvalue is one and hence the uncontrollable mode is stable.

As has been shown in Section~\ref{Sec:3}, system ($\widetilde{A},\widetilde{B}$) in \eqref{SEq:LinearTransformationS} is stabilizable. Note that the first column in $\widetilde{B}_{k}$, \emph{i.e.}, $\widetilde P_1$, is equal to $\widetilde{B}$, which indicates the stabilizability of system ($\widetilde{A},\widetilde P_1$). Then it is easy to observe that system ($\widetilde{A},\widetilde{B}_k$) is also stabilizable. According to Lemmas~\ref{Lemma:InvarianceLinear} and \ref{Lemma:InvarianceStateFeedback}, we can conclude that the mixed traffic system with multiple autonomous vehicles given by~\eqref{SEq:MultipleAVModel} is stabilizable.
\end{IEEEproof}

The specific expression of the uncontrollable mode is
\begin{equation*}  
\tilde {x}_{11}=\frac{\sum_{i\in S_{\text{AV}}}{\left(s_{i}(t)-s_{i,\text{c}}^{*}\right)}+\sum_{i\in \{
		1,2,\ldots ,n\} \backslash S_{\text{AV}}}{\left(s_{i}(t)-s^{*}\right)}}{\sqrt{n}},
\end{equation*}
which remains unchanged during the system evolution. The physical interpretation is the same as the case where there is one single autonomous vehicle: the sum of each vehicle's spacing should remain constant due to the ring road structure.

\begin{figure*}[t]
	\centering
	\setlength{\abovecaptionskip}{0em}
	\setlength{\belowcaptionskip}{0em}
	\newcommand{\fighspace}{\hspace{5mm}}
	\subfigure[]
	{ \label{Fig:UnstableWave}
		\includegraphics[scale=0.43]{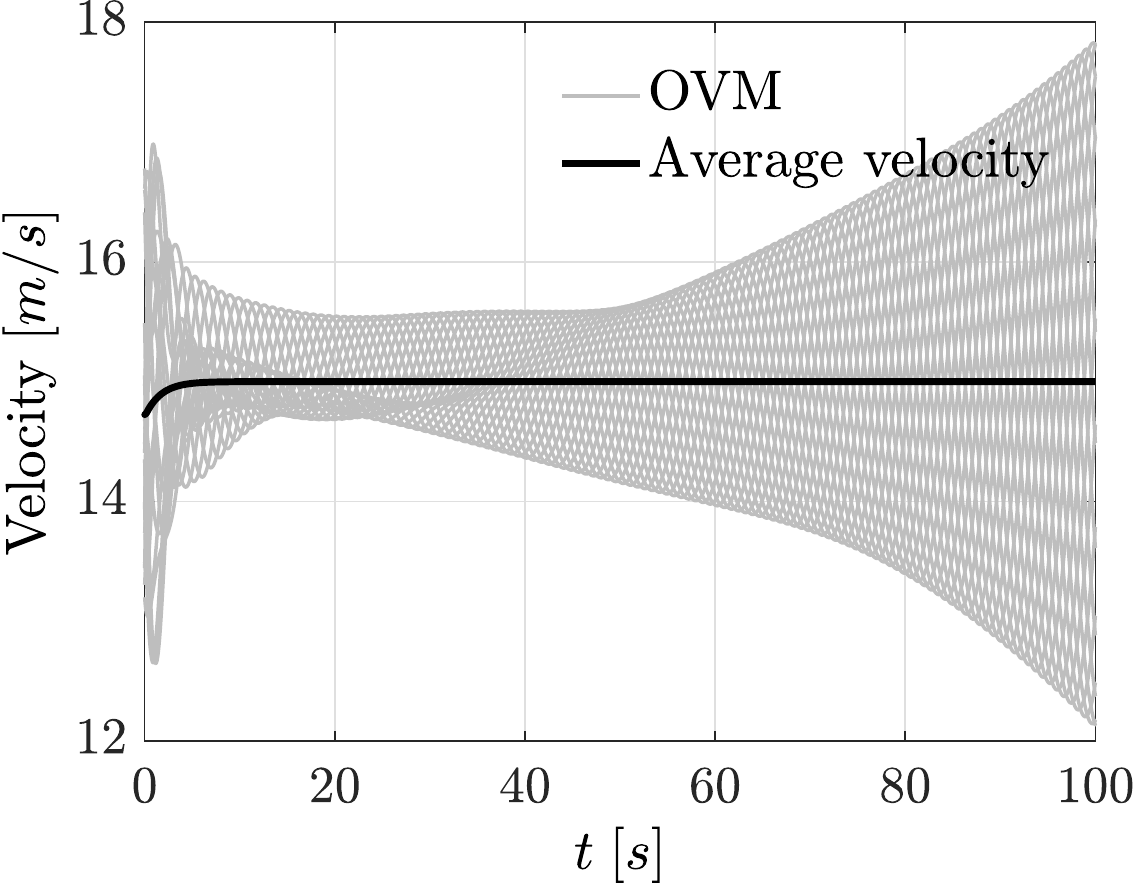}
	}
	\fighspace
	\subfigure[]
	{ \label{Fig:StabilizingFlowa}
		\includegraphics[scale=0.43]{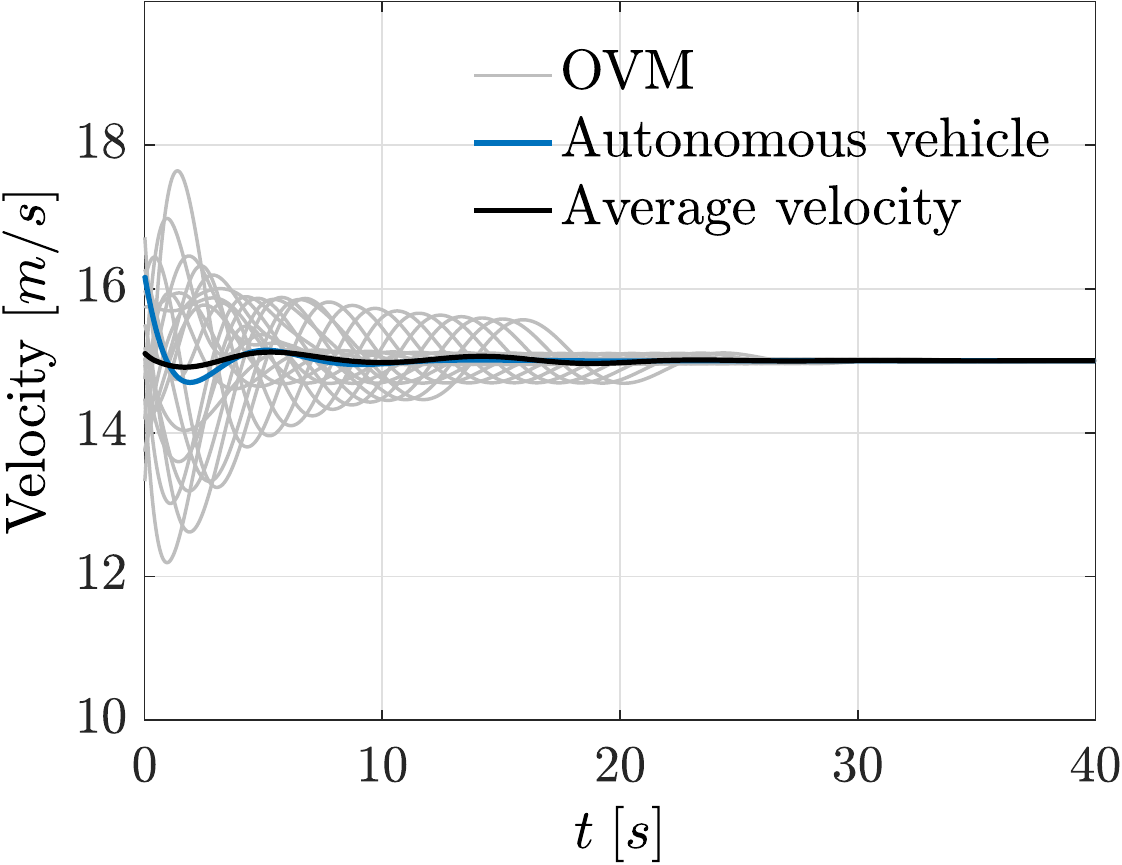}
	} 
	\fighspace
	\subfigure[]
	{ \label{Fig:StabilizingFlowb}
		\includegraphics[scale=0.43]{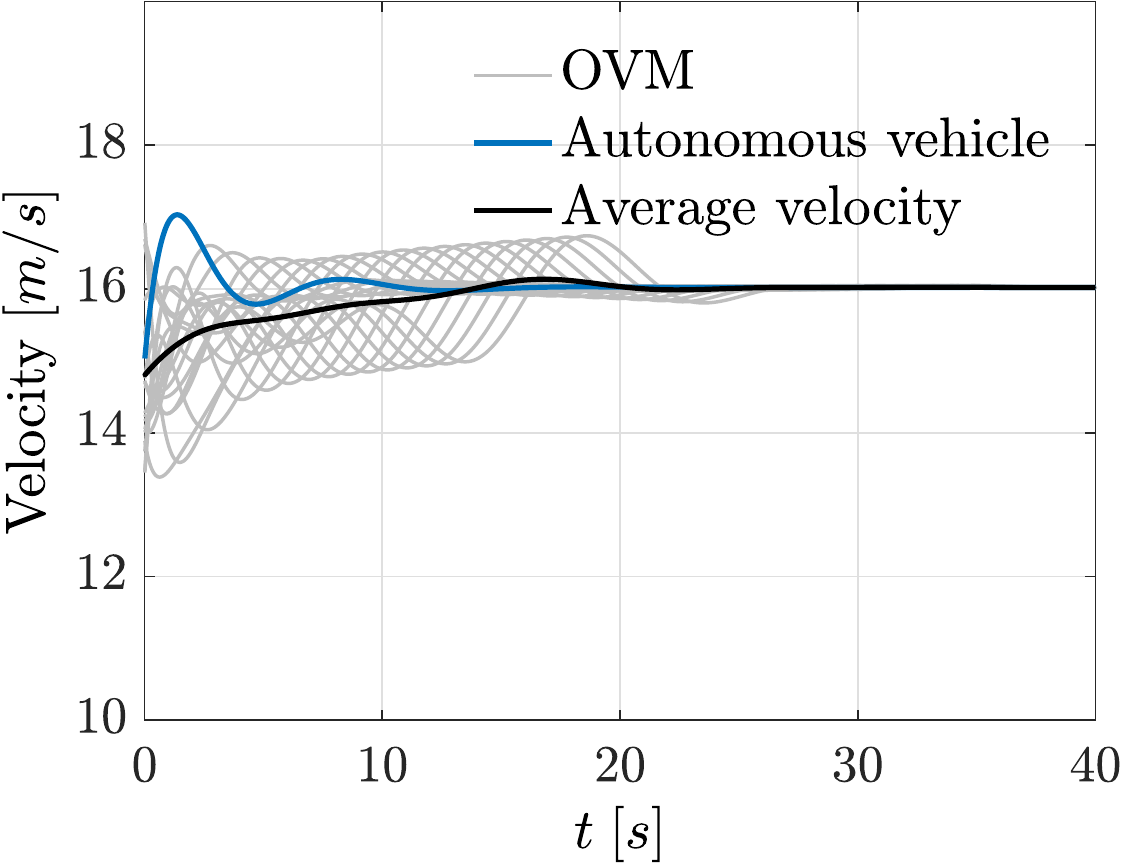}
	}
	\caption{Stabilizing traffic flow and increasing traffic speed. (a) The traffic system with human-driven vehicles only is unstable when $\alpha=0.6,\beta=0.9$ in the OVM model. (b) The mixed traffic system becomes stable after introducing an autonomous vehicle with an appropriate control strategy. (c) The traffic flow can be guided to a higher stable velocity (6\% improvement) via controlling the autonomous vehicle.}
	\label{Fig:StabilizingFlow}
\end{figure*}

\subsection{Optimal Control and Reachability Analysis}

Our controller formulation proposed in Section \ref{Sec:OptimalControl} can be easily applied to the mixed traffic system with multiple autonomous vehicles in \eqref{SEq:MultipleAVModel}. Define the feedback controller as $u(t)=-Kx(t)$, where
\begin{equation}
\label{Eq:FeedbackGainMultiple}
K=\begin{bmatrix}
K_1^\tr,K_2^\tr,\ldots,K_k^\tr
\end{bmatrix}^\tr,
\end{equation}
with $K_r \in \mathbb{R}^{1\times 2n}$ denoting the feedback gain for the autonomous vehicle indexed as $i_r$, \ie, $u_{i_r}(t)=-K_rx(t)$. Then following the process in Section \ref{Sec:OptimalControl}, we can obtain an optimal control strategy. 
%
%
We remark that the specific feedback gain $K_r$ for each autonomous vehicle may differ from each other depending on their positions in the traffic system. The resulting controller is a cooperation strategy for multiple AVs to achieve an optimal performance of the entire traffic system.

The reachability analysis is as follows.


\begin{theorem} 
	Consider the mixed traffic system in a ring road with multiple AVs in~\eqref{SEq:MultipleAVModel}. Suppose a stabilizing feedback gain \eqref{Eq:FeedbackGainMultiple} is found by \eqref{SEq:H2problemTraffic} and the coefficient matrix of~\eqref{Eq:ReachaiblityEquationMultiple} is non-singular. Then the traffic system can be stabilized at velocity $v^*$, if and only if the sum of the desired spacing $s_{i,c}^*$ ($i\in S_{\text{AV}}$) of each AV satisfies
	\begin{equation} \label{Eq:ReachabilityConditionMultiple}
	\sum_{i\in S_{\text{AV}}}s_{i,\text{c}}^*=L-(n-k)s^*,
	\end{equation}
	with $s^*$ given by the equilibrium equation \eqref{Eq:EquilibriumEquation}.
\end{theorem}

\begin{IEEEproof}
	Denote ($r=1,2,\ldots,n$)
	$$
	K_r= \begin{bmatrix} k^{(r)}_{1,1},k^{(r)}_{1,2},k^{(r)}_{2,1},k^{(r)}_{2,2,}, \ldots, k^{(r)}_{n,1},k^{(r)}_{n,2} \end{bmatrix}.
	$$
	Similar to the reachability analysis in Section \ref{Sec:4}, the final state of stable system in~\eqref{SEq:MultipleAVModel} can be obtained via
	\begin{equation*}
	\begin{cases}
	\tilde{v}_{i}(t_{\text{f}})=v_{\text{e}}, & i\in \{1,2,\ldots ,n \}, \\
	\tilde{s}_{i}(t_{\text{f}})=\displaystyle \frac{\alpha _{2}-\alpha _{3}}{\alpha _{1}}v_{\text{e}},
	& i\in \{ 1,2,\ldots ,n\} \backslash S_{\text{AV}},\\
	\tilde{s}_{i}(t_{\text{f}})=s_{\text{e},i}, & i\in S_{\text{AV}}.\\
	\end{cases}
	\end{equation*}
	where $v_\text{e}$ and $s_{\text{e},i},\; i\in S_{\text{AV}}$ should satisfy
	\begin{equation} \label{Eq:ReachaiblityEquationMultiple}
	\begin{cases}
	\eta^{(1)} v_e + \sum_{i\in S_{\text{AV}}} k^{(1)}_{i,1} s_{\text{e},i} = 0,    \\
	\quad \quad \quad \vdots    \\
	\eta^{(k)} v_e + \sum_{i\in S_{\text{AV}}} k^{(k)}_{i,1} s_{\text{e},i} = 0,    \\
	(n-k)\frac{\alpha_2-\alpha_3}{\alpha_1}  v_\text{e}  + \sum_{i\in S_{\text{AV}}}s_{\text{e},i} =L_k,
	\end{cases}
	\end{equation}
	with ($r=1,2,\ldots,n$)
	$$
	\begin{aligned}
	&\eta^{(r)} =\frac{\alpha_2-\alpha_3}{\alpha_1} \displaystyle \sum_{i\in \{1,2,\ldots,n\} \backslash S_{\text{AV}}} k^{(r)}_{i,1} + \displaystyle \sum_{i \in \{1,2,\ldots,n\}} k^{(r)}_{i,2} ,\\
	&L_k=L-(n-k)s^*-\sum_{i\in S_{\text{AV}}}s_{i,\text{c}}^*.
	\end{aligned}
	$$
	
	The condition that the traffic flow reaches the desired equilibrium velocity $v^*$ is equivalent to $v_\text{e}=0,\;s_{\text{e},i}=0, i\in S_{\text{AV}}$. The solution to \eqref{Eq:ReachaiblityEquationMultiple} is all zeros, \ie, $v_\text{e}=0$ and $s_{\text{e},i}=0,\; i\in S_{\text{AV}}$, if and only if the constant vector is zero, \emph{i.e.}, $L_k=0$, which leads to \eqref{Eq:ReachabilityConditionMultiple}.
\end{IEEEproof}

Since the spacing of each autonomous vehicle cannot be negative, \textit{i.e.}, $s_{i,c}^*>0$ ($i\in S_{\text{AV}}$), we have
$
\sum_{i\in S_{\text{AV}}}{s_{i,c}^{*}}=L-(n-k)s^{*}>0.
$
In this case, the maximum spacing for each HDV in the equilibrium can be increased to
 \begin{equation*}  
	s^*_{\max,k} < \frac{L}{n-k},
	\end{equation*}
which sets up a new maximum equilibrium traffic velocity $(v^*)_{\max,k}$. 
 \begin{corollary}
	There exists a reachable range for the traffic velocity in the mixed traffic system with $k$ autonomous vehicles given by \eqref{SEq:MultipleAVModel}, \emph{i.e.},
	$0\leqslant v^* < v^*_{\max,k},$
	where
	\begin{equation*}
	   F\left(\frac{L}{n-k},0,v^*_{\max,k}\right)=0.
	\end{equation*}
\end{corollary}

 This result generalizes the statement in Corollary \ref{Corollary:Velocity}, and it can be observed clearly that a larger proportion of autonomous vehicles leads to a higher reachable traffic velocity.

\section{Numerical Experiments}
\label{Sec:6}
Our theoretical results are obtained using a linearized model of the mixed traffic system. We evaluate their effectiveness in the presence of nonlinearities arising in the car-following dynamics. In this section, we conduct multiple simulation experiments to validate our results based on a realistic nonlinear OVM model, as shown in \eqref{Eq:OVM}. 
All the experiments are carried out in MATLAB.

\subsection{Experimental Setup}

Similarly to \cite{jin2017optimal}, we set the parameters in the OVM model \eqref{Eq:OVM} as follows: $\alpha=0.6$, $\beta=0.9$, $s_{\mathrm{go}}=35$, $v_{\max}=30$, $s_{\mathrm{st}}=5$. One can verify that this set of values violates the stability condition \eqref{Eq:StabilityConditionforOVM}, which means that a ring-road traffic system with such HDVs only is unstable and stop-and-go waves may happen in case of any perturbations.

For the parameters in the performance output \eqref{Eq:Output}, we choose $\gamma_s=0.03$, $\gamma_v=0.15$, $\gamma_u=1$. Based on the approach in Section \ref{Sec:4}, an optimal linear feedback gain $K$ for the autonomous vehicle is obtained using Mosek~\cite{mosek2010mosek}.
To avoid crashes, we also assume that all the vehicles are equipped with a standard automatic emergency braking system
\begin{equation*}
	\dot{v}(t)=a_{\min },\,\text{if}\,\frac{v_{i}^{2}(t)-v_{i-1}^{2}(t)}{2s_{i}(t)}\ge \vert a_{\min }\vert,
\end{equation*}
where the maximum deceleration rate of each vehicle is set to $a_{\min}=-5m/s^2$.

\subsection{Stabilizing Traffic Flow and Increasing Traffic Velocity}

Our first experiment aims to show that a typical nonlinear mixed traffic flow can be stabilized by one single autonomous vehicle, as proved in Theorem~\ref{Theorem:stabilizability} after linearization. We consider the case where $n=20$ and $L=400$. Each vehicle has a weak perturbation around its equilibrium state at initial time, in the sense that the position and the velocity of the $i$-th vehicle are $iL/n+\delta s$, and $v_{\text{ini}}+\delta v$, where $v_{\text{ini}}$ is the equilibrium velocity corresponding to the equilibrium spacing $L/n$, $\delta s\sim U[-4,4]$ and $\delta v\sim U[-2,2]$ with $U[a,b]$ denoting a uniform distribution between $a$ and $b$.

 When all the vehicles are under human control, it is clearly observed in Fig.\ref{Fig:UnstableWave} that the initial perturbations inside the traffic flow are amplified gradually, and each vehicle's velocity keeps fluctuating, finally inducing a stop-and-go wave. In contrast, if there is one autonomous vehicle using the proposed control method, the traffic flow can be stabilized to the original average velocity $15m/s$ within a short time (Fig.\ref{Fig:StabilizingFlowa}). Moreover, by adjusting the equilibrium velocity $v^*$ (as stated in Section~\ref{Sec:Reachability}) and the corresponding equilibrium spacing $s^*$ and $s_c^*$ according to \eqref{Eq:EquilibriumEquation} and \eqref{Eq:ReachabilityCondition} respectively, the AV can steer the entire traffic flow towards a higher velocity, from $15m/s$ to $16m/s$; see Fig.\ref{Fig:StabilizingFlowb}. In this case we observed 6\% improvement of traffic velocity when there exist only 5\% AVs (one out of 20) in the mixed traffic system.

\begin{figure}[t]
	\centering
	\setlength{\abovecaptionskip}{0em}
	\setlength{\belowcaptionskip}{0em}
	\subfigure[]
	{ \label{Fig:ControlEnergy_Comparison}
		\includegraphics[scale=0.4]{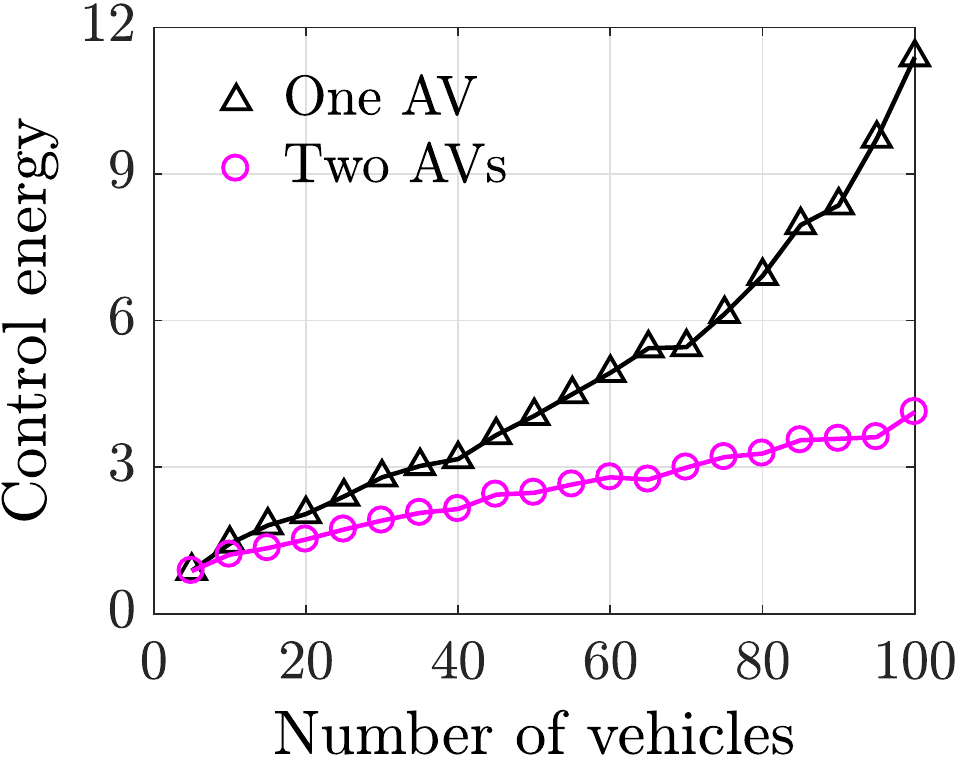}
	}
	\subfigure[]
	{ \label{Fig:TimeToStabilize_Comparison}
		\includegraphics[scale=0.4]{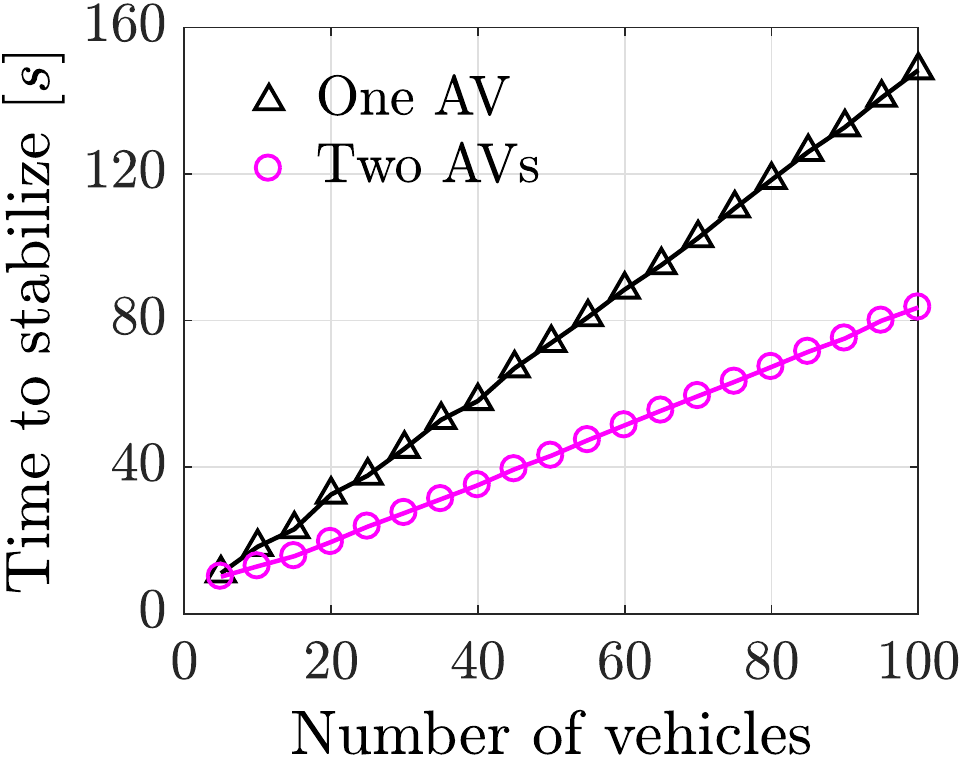}
	}
	\caption{Simulation results at different system scales. We ran 2000 random simulations for each value of $n$. The parameters are as follows: $\gamma_s=0.03$, $\gamma_v=0.15$, $\gamma_u=1$. (a) The control energy $\int_0^{\infty} u^\tr u dt$ needed to stabilize the traffic flow for each autonomous vehicle. (b) The time required to stabilize the traffic system. }
	\label{Fig:SystemSize}
\end{figure}

\subsection{Smoothing Traffic Flow via Multiple Autonomous Vehicles}

Fig.\ref{Fig:TrafficWave} and Fig.\ref{Fig:StabilizingFlow} have demonstrated the ability of a single autonomous vehicle to smooth the traffic flow where there exist certain perturbations. We proceed to conduct numerical experiments for the scenario where the traffic system has two autonomous vehicles. We consider different values of $n$ and let $L=20n$. The experimental setup at initial time is the same as that in the previous experiment. The results are shown in Fig.\ref{Fig:SystemSize}. It is clear that both the settling time and the control energy of each autonomous vehicle decrease by a factor of two approximately, when there are two autonomous vehicles in the traffic system uniformly. Based on the results, we may estimate the market penetration rate of autonomous vehicles to control traffic flow effectively when adopting the optimal control strategy. In the scenario of Fig.\ref{Fig:SystemSize}, if one wants to reject the influence of the perturbation on traffic flow within 30 seconds, a single autonomous vehicle can control the traffic flow consisting of around 20 HDVs. This number 
agrees with the results from real-world experiments~\cite{stern2018dissipation}.

\subsection{Dampening Traffic Waves and Comparison with Existing Strategies}

 As our last experiment, we consider a scenario with the presence of infrastructure bottlenecks or lane changing~\cite{stern2018dissipation}, where one vehicle has a rapid deceleration representing a strong perturbation. In the beginning, the traffic flow is at the equilibrium state with the velocity $15 m/s$. And then at $t = 20s$, the $i$-th vehicle decelerates to $5 m/s$ in two seconds. We observe that if all the vehicles are human-driven, the perturbation may grow stronger during the propagation process (Fig.\ref{Fig:StrongPerturbationHDV_2}), while the autonomous vehicle with an optimal control strategy can respond actively to attenuate the perturbation and stabilize the traffic flow (Fig.\ref{Fig:StrongPerturbationAV_2}). Here we only show the case where the 6th vehicle is under the strong perturbation. Indeed, the experiment results confirm that our strategy allows one autonomous vehicle to dampen strong traffic waves wherever they come from.

\begin{figure}[t]
	\centering
	\setlength{\abovecaptionskip}{0em}
	\setlength{\belowcaptionskip}{0em}
	\newcommand{\fighspace}{\hspace{1mm}}
	\subfigure[]
	{ \label{Fig:StrongPerturbationHDV_2}
		\includegraphics[scale=0.42]{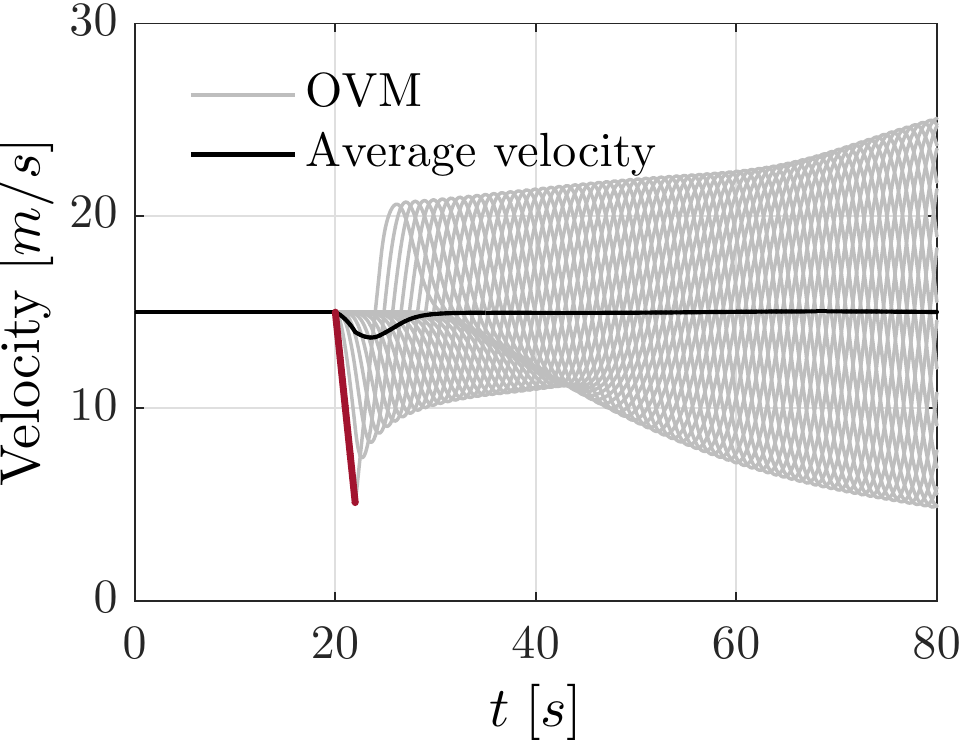}
		\includegraphics[scale=0.42]{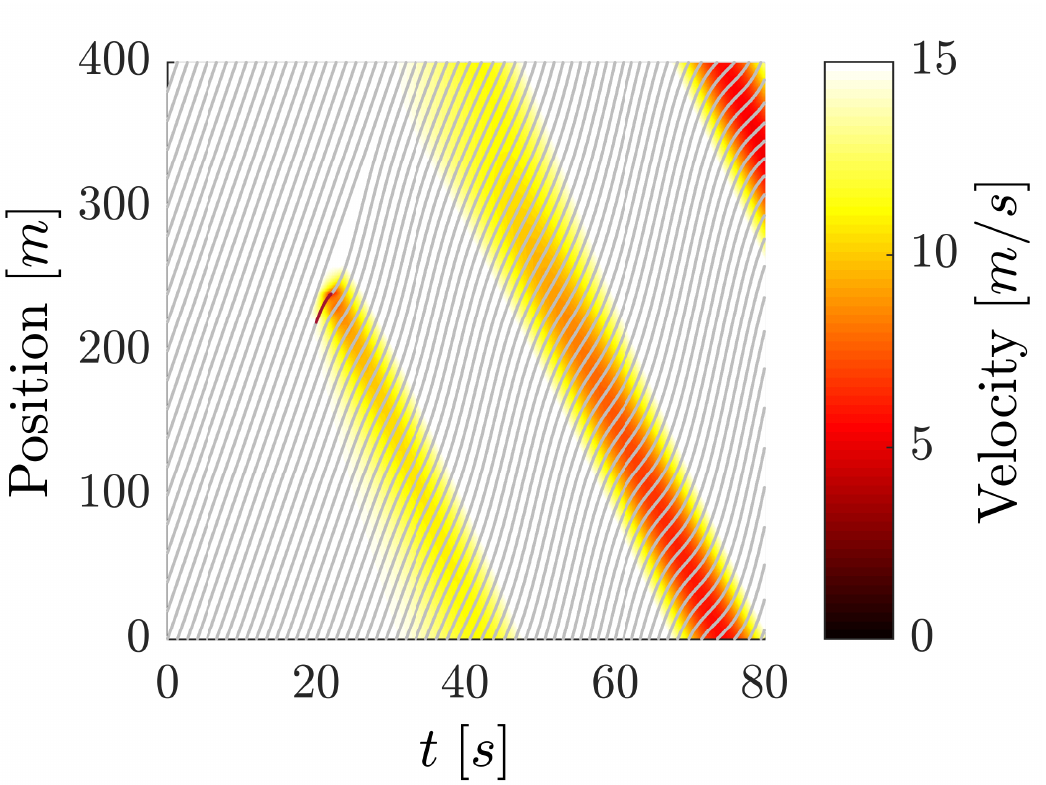}
	}
	\subfigure[]
	{ \label{Fig:StrongPerturbationAV_2}
		\includegraphics[scale=0.42]{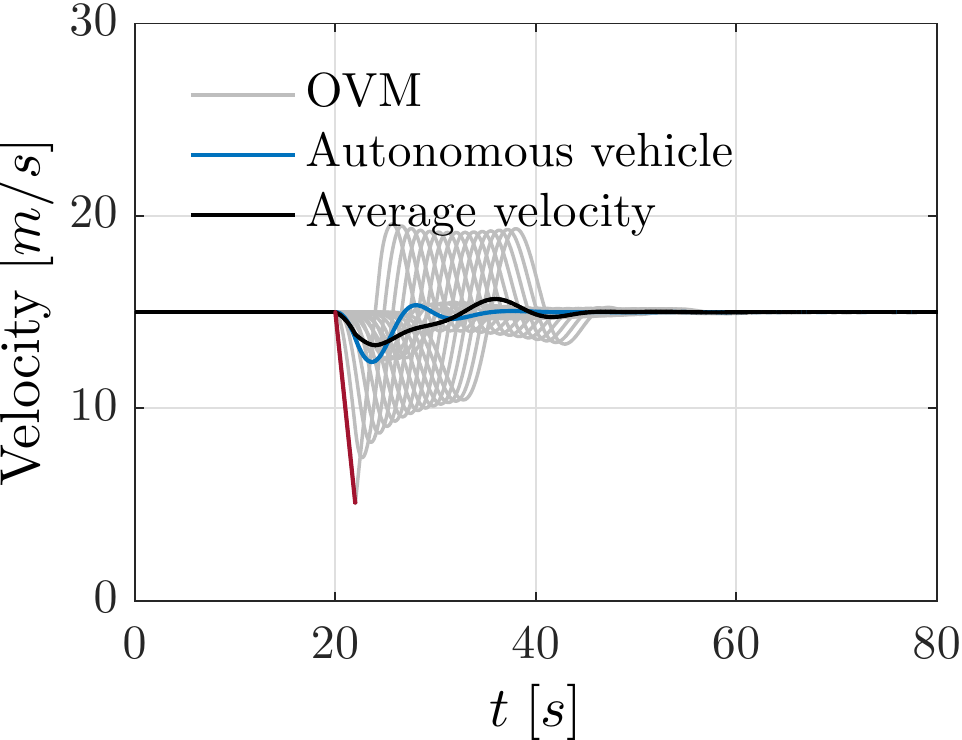}
		\includegraphics[scale=0.42]{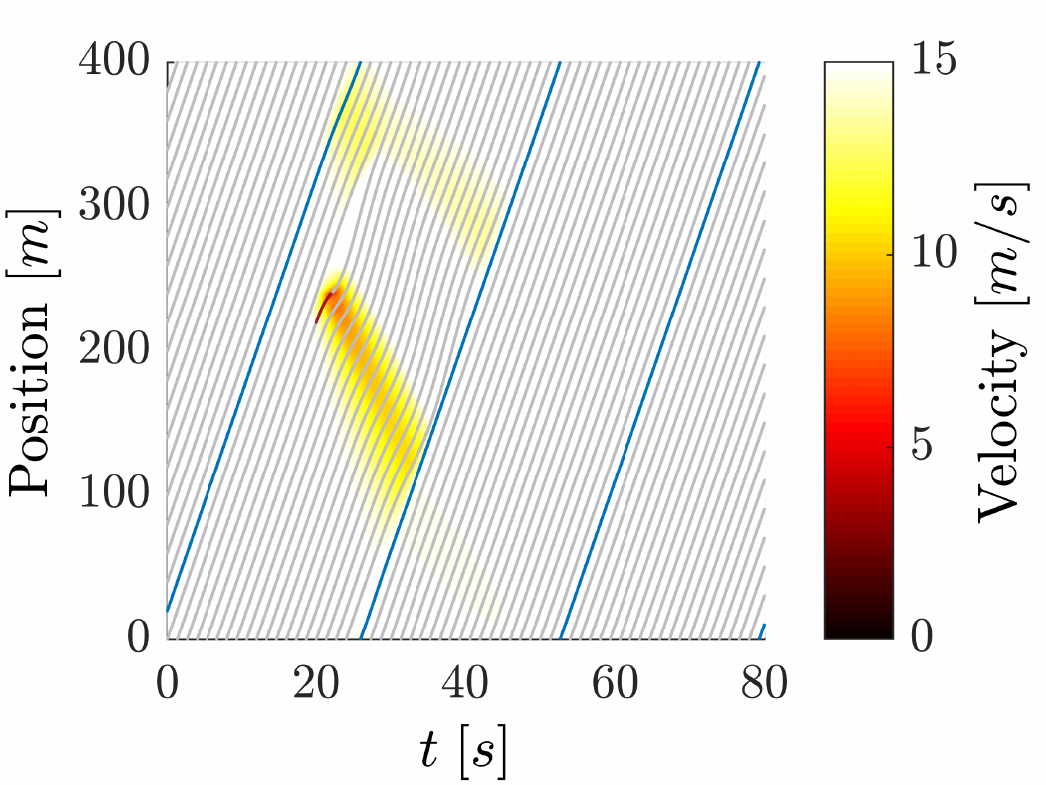}
	}
	\caption{Numerical results for the scenario with a rapid and strong perturbation in the $6$th vehicle. (a) The traffic system consists of HDVs only. (b) The mixed traffic system has an autonomous vehicle that adopts the optimal control strategy. In each panel, the right figure shows the vehicles' trajectories, where the red zone represents the traffic wave; the left figure shows the vehicles' velocities, where the red line denotes the perturbation and the black line is the average velocity of all vehicles. }
	\label{Fig:StrongPerturbation}
\end{figure}

\begin{figure}[t]
	\centering
	\setlength{\abovecaptionskip}{0em}
	\setlength{\belowcaptionskip}{0em}
	\includegraphics[scale=0.48]{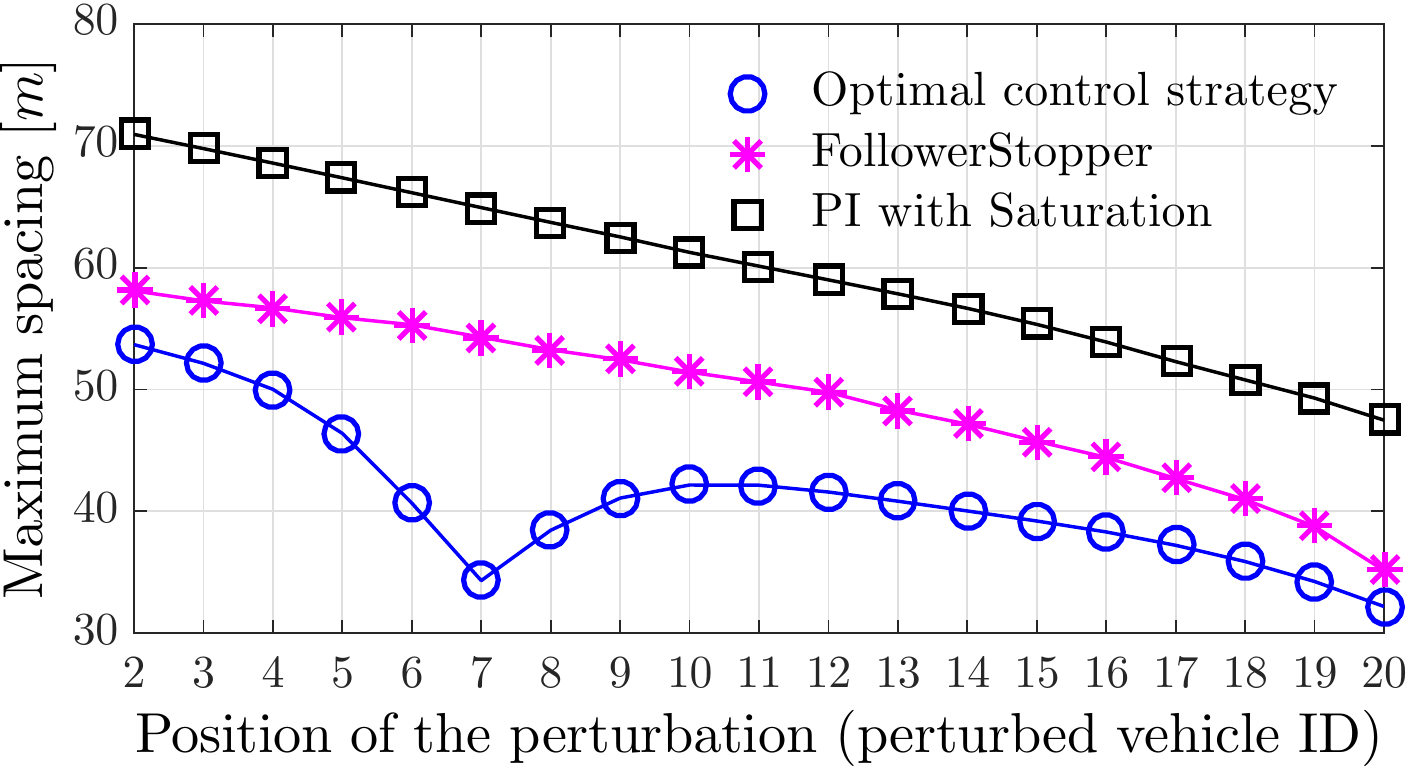}
	\caption{ Comparison between three strategies of the maximum spacing of the autonomous vehicle, \ie, $\max_{t}{s_1(t)} $, during the experiment with a rapid and strong perturbation.}
		\label{SFig:Spacing_Comparison}
		\vspace{-3mm}
\end{figure}

Next, we compare our proposed strategy with existing ones. As shown in Fig.\ref{Fig:TrafficWave}, our strategy allows the autonomous vehicle to mitigate undesired perturbations in an active way instead of responding passively as CACC-type controllers. Here, we proceed to make comparisons with two heuristic controllers that aim at dampening traffic waves:
FollowerStopper and PI with Saturation \cite{stern2018dissipation}. Note that the FollowerStopper and PI with Saturation use a command velocity $v_\mathrm{cmd}$ as the control input, and we add a proportional controller, \ie, $u(t)=k_p (v_\mathrm{cmd} (t)-v(t))$, with $k_p=0.6$, to serve as a lower-level controller. In all tested scenarios, the perturbation was successfully dampened using the three methods. However, since FollowerStopper and PI with Saturation are essentially slow-in fast-out strategies, they tend to leave a long gap from the preceding vehicle, which may cause vehicles from adjacent lanes to cut in. In contrast, our optimal strategy avoids this problem and keeps the spacing within a moderate range. As shown in Fig.\ref{SFig:Spacing_Comparison}, this finding holds irrespectively of the position of the perturbation, and hence confirms the advantage of our strategy. 

\begin{figure}[t]
	\centering
	\setlength{\abovecaptionskip}{0em}
	\setlength{\belowcaptionskip}{0em}
	\includegraphics[scale=0.46]{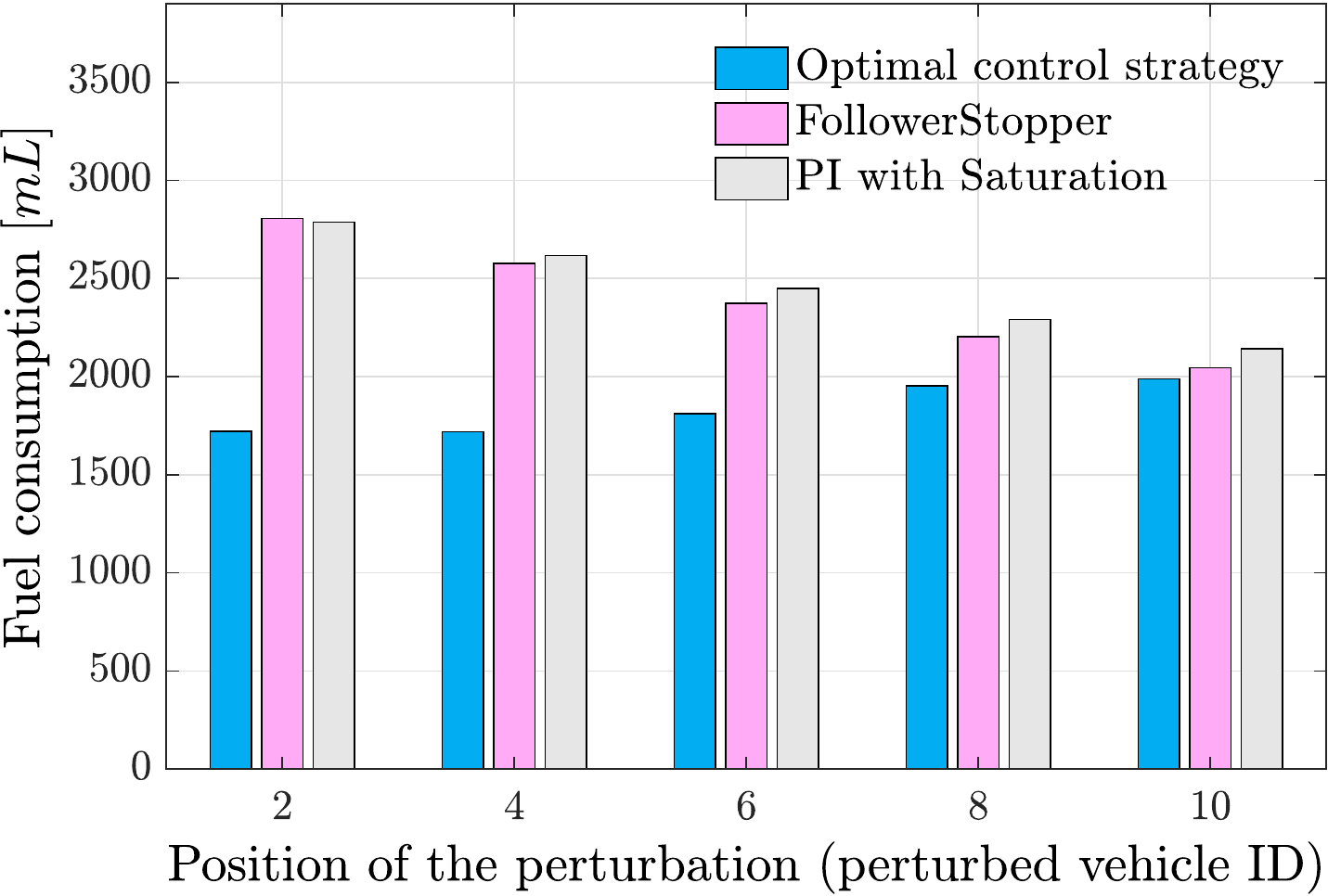}
	\caption{Comparison between three strategies of the total fuel consumption of the entire traffic flow, \ie, $FC$ in \eqref{Eq:Total_FC_model}, during the experiment with a rapid and strong perturbation.}
	\label{SFig:FuelConsumption_Comparison}
\end{figure}

Finally, we consider the comparison of fuel consumption for the three control methods. It has been demonstrated in \cite{stern2018dissipation} that applying FollowerStopper or PI with Saturation to one autonomous vehicle can result in a $40\%$ reduction of the total fuel consumption of the entire traffic flow. Here we are interested in whether our strategy can achieve further improvement. An instantaneous fuel consumption model in~\cite{bowyer1985guide} is utilized to estimate the fuel consumption rate $f_i$ ($mL/s$) of the $i$-th vehicle, which is given by
\begin{equation*}
	f_i = \begin{cases}
	0.444+0.090 R_i v_i + [0.054 a_i^2 v_i] _{a_i>0},& \text{if}\; R_i>0,\\
	0.444, & \text{if} \; R_i \le 0,
	\end{cases}
\end{equation*}
where $R_i = 0.333+0.00108 v_i^2 + 1.200 a_i$. Then, the total fuel consumption of the entire traffic flow $FC$ is calculated as
\begin{equation} \label{Eq:Total_FC_model}
FC = \sum_{i=1}^{N}\left({\int_{t=0}^{t_f}f_i dt}\right),
\end{equation}
where $t_f=100s$ denotes the end time of the simulation. The simulation setup is the same as that in Fig.\ref{Fig:StrongPerturbation}, and the results are shown in Fig.\ref{SFig:FuelConsumption_Comparison}. It is observed that our proposed strategy achieves evidently lower fuel consumption than FollowerStopper and PI with Saturation when the perturbation happens within the range from the $1$st to the $10$th vehicle. This result validates the great potential of our strategy in improving fuel economy. We note that, when the perturbation happens within the range from the $11$th to the $20$th vehicle, which is ahead of the autonomous vehicle in a small distance, all of the three strategies require the autonomous vehicle to brake hard to guarantee safety. In these cases, the three control strategies have similar performance in terms of fuel consumption.

\section{Conclusion}
\label{Sec:7}

Unlike the traditional control methods that regulate traffic flow externally at fixed positions, autonomous vehicles can be used as {mobile} actuators to control traffic flow internally.
In this paper, we have introduced a comprehensive theoretical analysis to address the potential of autonomous vehicles on smoothing mixed traffic flow. Specifically, we have analyzed the controllability, stabilizability, and reachability of mixed traffic systems. Also, an optimal control strategy has been introduced to {actively} smooth mixed traffic flow.

A few other topics are worth further investigations. First, we have assumed autonomous vehicles have access to the global traffic state, \textit{i.e.}, the information of all other human-driven vehicles. Due to the limit of communication ranges, autonomous vehicles may be only able to obtain the information of its neighboring vehicles. It is interesting to design a localized optimal controller, and this leads to the notion of structured controller synthesis~\cite{zheng2018scalable,jovanovic2016controller}. Second, we have assumed homogeneous dynamics for human-driven vehicles, and potential time delays are ignored. One interesting direction is to consider heterogeneity and time delay in controlling mixed traffic systems. We note that some recent work has considered the effect of heterogeneity and time delays at the level of platoon control~\cite{di2015distributed, gao2016robust, jin2017optimal}, which may offer some insights for controller design in mixed traffic systems. Finally, our current analysis {focuses} on the single-lane ring road setting, and it would be interesting to extend our analysis to the scenarios with multiple lanes and lane-changing behavior.


%

\appendices
\section{Diagonalization of Block Circulant Matrices}
\label{Appendix:Circulant}

Define $\omega=e^{\frac{2\pi j}{n}}$, where $j=\sqrt{-1}$ denotes the imaginary unit, and the Fourier matrix $F_n$ is defined as~\cite{marshall2004formations,olson2014circulant}

\begin{equation}  \label{Eq:FourierMatrix}
F_n^*=\frac{1}{\sqrt{n}}\begin{bmatrix} 1 & 1 & 1 &\ldots & 1 \\
1 & \omega & \omega^2 & \ldots & \omega^{n-1}  \\
1 & \omega^2 & \omega^4 & \ldots & \omega^{2(n-1)}  \\
\vdots & \vdots & \vdots &  & \vdots\\
1 & \omega^{n-1} & \omega^{2(n-1)} & \ldots & \omega^{(n-1)(n-1)}  \\
\end{bmatrix},
\end{equation}
where $F_{n}^*$ denotes the conjugate transpose matrix of $F_{n}$. By the definition of Fourier matrix, we know that $F_n$ and $F_n^*$ are symmetric, \textit{i.e.}, $F_{n}^{*}=(F_{n}^{*})^\tr$, $F_{n}=F_{n}^\tr$, and that $F_{n}$ is a unitary matrix, \textit{i.e.}, $F_{n}F_{n}^{*}=I_{n}$, where $I_n$ denotes the $n \times n$ identity matrix.

Given $M_1,M_2,\ldots,M_n\in \mathbb{R}^{m\times m}$, a block circulant matrix is of the following form
\begin{equation} \label{Eq:BlockCirculantMatrix}
	M = \begin{bmatrix}
	A_1 & A_2 & \ldots & A_n\\
	A_n & A_1 & \ldots & A_{n-1}\\
	\vdots & \vdots &  & \vdots\\
	A_2 & A_3 & \ldots & A_1
	\end{bmatrix}\in \mathbb{R}^{mn\times mn}.
\end{equation}
As shown in~\cite{marshall2004formations,olson2014circulant}, a block circulant matrix given by~\eqref{Eq:BlockCirculantMatrix} can be diagonalized as follows
\begin{equation*}
	\text{diag}(D_{1},D_{2},\ldots
	,D_{n}) = (F_{n}^{*}\otimes I_{m})^{-1} \hat A (F_{n}^{*}\otimes I_{m}),
\end{equation*}
where
$$
\begin{bmatrix}
D_1\\D_2\\ \vdots \\D_n
\end{bmatrix}=
(\sqrt{n}F_n^*\otimes I_m)\begin{bmatrix}
M_1\\M_2\\ \vdots \\M_n
\end{bmatrix}.
$$

\section{Proof of the Stability Condition \eqref{Eq:StabilityCondition}} \label{appendix:stability}
When all the vehicles are controlled by human drivers, the system dynamics is given by
\begin{equation}
	\dot x = \hat{A} x,
\end{equation}
where $\hat{A}$ is the same as \eqref{Eq:Ahat}. To analyze the stability of matrix $\hat{A}$, it is necessary and sufficient to study its eigenvalues' distribution. Since $\hat{A}$ is a block circulant matrix, it can be diagonalized to simplify the eigenvalue calculation. As shown in Appendix \ref{Appendix:Circulant}, $\hat{A}$ can be diagonalized into
 \begin{equation*}
	\hat A=(F_n^* \otimes I_2) \cdot \text{diag}(D_1,D_2,\ldots ,D_n)\cdot (F_n \otimes I_2),
	\end{equation*}
where $D_{i}=A_{1}+A_{2}\omega^{(n-1)(i-1)}, i=1, 2, \ldots, n$. Since
\begin{equation}
	\begin{aligned}
	\text{det}(\lambda I-A)&=\text{det}(\lambda I-\text{diag}(D_{1},D_{2},\ldots, D_{n}))\\
&	=\prod_{i=1}^{n}{\text{det}(\lambda I-D_{i})},
	\end{aligned}
\end{equation}
then the eigenvalues $\lambda$ of $\hat{A}$ can be calculated by
	\begin{equation} \label{AEq:EigenvalueEquation}
	\lambda^2+\left(\alpha_2-\alpha_3 \omega^{(n-1)(i-1)}\right)\lambda+\alpha_1\left(1-\omega^{(n-1)(i-1)}\right)
	=0,
	\end{equation}
where $i=1,2,\ldots,n$. Note that \eqref{AEq:EigenvalueEquation} is a second-order complex equation, which makes it non-trivial to get the analytical roots. Instead, we transform this equation to continue our analysis. Substituting the expression of $\omega $ into~\eqref{AEq:EigenvalueEquation} leads to
	\begin{equation} \label{Eq:EigenvalueH}
	\displaystyle e^{\frac{i-1}{n}\cdot 2\pi j}=\frac{\alpha _{1}+\alpha _{3}\lambda }{\alpha _{1}+\alpha _{2}\lambda +\lambda ^{2}}=H(\lambda ),\quad i=1,
	2,\ldots ,n,
	\end{equation}
which means that the eigenvalues of $\hat{A}$ correspond to the solutions of~\eqref{Eq:EigenvalueH}. Note that $e^{\frac{i-1}{n}\cdot 2\pi j}$ is the $i$-th complex root of $z^{n}=1$, indicating that for all the eigenvalues $\lambda $ of ${A}$, the values of $H(\lambda )$ constitute $n$ unit roots. As $n$ changes, $H(\lambda )$ corresponds to different unit roots. Therefore, if all the roots of $ \vert H(\lambda )\vert =1$ have negative real parts, then the solutions of equation~\eqref{Eq:EigenvalueH}, \emph{i.e.}, the eigenvalues of matrix $\hat{A}$, have negative real parts. We conclude that the condition that all the roots of $\vert H(\lambda )\vert =1$ have negative real parts is sufficient to guarantee that $\hat{A}$ is stable. Note that this condition becomes
sufficient and necessary for the case where the system is stable for any $n$. This is because that $H(\lambda )$ can be any unit root $e^{\theta j}$, for $\theta \in[0,2\pi) $.

The rest of analysis is the same as~\cite{cui2017stabilizing}. Since all rational functions are meromorphic, $H(\lambda)$ is a meromorphic function. Because $\alpha _{1}$ and $\alpha _{2} $ are positive real numbers, the poles of $H(\lambda )$ are in the left half plane, indicating that $H(\lambda )$ is holomorphic in the right half plane. Meanwhile, $\vert H(\lambda )\vert \to 0$ when $\text{Re}(\lambda )\to \infty $. According to Maximum Modulus Principle~\cite{burckel1980introduction}, the extreme value of $\vert H(\lambda )\vert $ in the right half plane can only be obtained on the imaginary axis.  To avoid eigenvalues with positive real parts, $\vert H(\lambda )\vert $ should not be more than $1$ on the imaginary axis. Therefore, that the roots of $\vert H(\lambda )\vert =1 $ have negative real part is equivalent to
$
\vert H(jv)\vert \le 1,\  \forall v\in \mathbb{R}.
$ This inequality leads to the stability criterion
$	\alpha _{2}^{2}-\alpha _{3}^{2}-2\alpha _{1}\ge 0.$
%

\section*{Acknowledgment}

The authors thank Prof.~Antonis Papachristodoulou, Prof.~Jianqiang Wang, Dr.~Qing Xu, Mr.~Licio Romao, Mr.~Suhao Yan and Mr.~Chaoyi Chen for comments and discussions. We would also like to thank Dr.~Ross Drummond at the University of Oxford for providing constructive feedback.

\ifCLASSOPTIONcaptionsoff
  \newpage
\fi



%

\bibliographystyle{IEEEtran}
\bibliography{IEEEabrv,mybibfile}

\end{document}